# CYCLE CARACTÉRISTIQUE SUR UNE PUISSANCE SYMÉTRIQUE D'UNE COURBE ET DÉTERMINANT DE LA COHOMOLOGIE ÉTALE

Fabrice Orgogozo[†] et Joël Riou[‡]

Résumé. En nous appuyant de manière essentielle sur le formalisme développé par Alexander Beilinson et Takeshi Saito, nous calculons le cycle caractéristique d'une puissance tensorielle symétrique externe d'un faisceau étale modéré sur une courbe. Ceci généralise un résultat de Gérard Laumon en caractéristique nulle et entraîne un résultat de locale acyclicité du morphisme d'Abel-Jacobi, dû à Pierre Deligne et motivé par son approche géométrique de la formule du produit pour le déterminant de la cohomologie (facteur epsilon).

Peu après le dépôt sur arXiv, Will Sawin nous a informé avoir obtenu une formule plus générale pour le cycle caractéristique que celle obtenue (indépendamment) dans ce texte.

Abstract. Relying on the formalism developed by Alexander Beilinson and Takeshi Saito, we compute the characteristic cycle of an external symmetric power of a tame étale sheaf on a curve. This generalizes a result of Gérard Laumon in characteristic 0 and leads to a result of local acyclicity of the Abel-Jacobi morphism, due to Pierre Deligne and motivated by his geometric approach to the product formula for the determinant of cohomology (epsilon factor).

Shortly after the submission on arXiv, Will Sawin kindly informed us that he had obtained a more general formula for the characteristic cycle than the one which we obtained (independently) in this text.

## Table des matières

## Introduction

**0.1**. Soient $Y$ un schéma algébrique lisse sur un corps algébriquement clos $k$, purement de dimension $n$, et un corps fini $\Lambda$ de caractéristique inversible dans $k$. Alexander Beilinson a associé à $\mathscr{G} \in \mathrm{D}^{\mathrm{b}}_{\mathrm{c}}(Y_{\text{ét}}, \Lambda)$ un fermé $SS(Y, \mathscr{G})$ du fibré cotangent $T^\star Y$ de $Y$, appelé support singulier de $\mathscr{G}$, purement de dimension $n$, qui contrôle en termes différentiels (une condition de transversalité) l'acyclicité des paires $(Y \to T, \mathscr{G})$, où $T$ est un $k$-schéma algébrique lisse ([Beilinson 2017, §1.2-3]). Takeshi Saito a raffiné cette construction en un *cycle* caractéristique $CC(Y, \mathscr{G}) \in Z^n(T^\star Y)$, dont le support est contenu dans $SS(Y, \mathscr{G})$, satisfaisant une formule de l'indice : si $Y$ est de plus projective, la caractéristique d'Euler-Poincaré de $\mathscr{G}$ est égale au degré d'intersection de ce cycle avec la section nulle du fibré cotangent ([Saito 2017b]).

[†]Institut de mathématiques de Jussieu—Paris rive gauche (IMJ-PRG), CNRS, 4 place Jussieu, 75005 Paris, France
Fabrice.Orgogozo+math@normalesup.org

[‡]Université Paris-Saclay, CNRS, Laboratoire de mathématiques d'Orsay, 91405, Orsay, France ; joel.riou@universite-paris-saclay.fr

Dans cet article, nous considérons de tels cycles lorsque $Y$ est le produit symétrique d'une courbe $X$ et $\mathscr{G}$ provient par l'opération « tenseurs symétriques » d'un faisceau *modéré* $\mathscr{F}$ sur $X$. Plus précisément, nous associons à $\mathscr{F}$ une série formelle

$$S_{\mathscr{F}} := \sum_{n\in\mathbf{N}} CC(X^{\langle n\rangle}, \mathscr{F}^{\langle n\rangle}) \in \prod_{n\in\mathbf{N}} Z^n(T^\star X^{\langle n\rangle}).$$

En nous appuyant de manière essentielle sur le formalisme développé par T. Saito ([Saito 2017a], [Saito 2021]) nous dégageons une structure multiplicative (partiellement définie) permettant des dévissages : si $0 \to \mathscr{F}' \to \mathscr{F} \to \mathscr{F}'' \to 0$ est une suite exacte de faisceaux sur $X$, dans le cas favorable, on a $S_{\mathscr{F}} = S_{\mathscr{F}'} \cdot S_{\mathscr{F}''}$. Nous obtenons ainsi une formule en termes des séries associées au faisceau constant *de rang* 1 et à des faisceaux gratte-ciels :

$$S_{\mathscr{F}} = (S_\Lambda)^r \cdot \prod_{s\in\Sigma} S_{\Lambda_s}{}^{-a_s},$$

où $r$ est le rang générique de $\mathscr{F}$ et les $a_s := r - \dim_\Lambda \mathscr{F}_s \in \mathbb{Z}$ sont les chutes du rang en les points singuliers $s \in S$ du faisceau modéré $\mathscr{F}$.
Pour expliciter ces formules, c'est-à-dire déterminer pour chaque $n$ le cycle $CC(X^{\langle n\rangle}, \mathscr{F}^{\langle n\rangle})$ (1.5.1), il apparaît essentiel de calculer la structure multiplicative partiellement définie

$$Z^n(T^\star X^{\langle n\rangle}) \times Z^m(T^\star X^{\langle m\rangle}) \dashrightarrow Z^{n+m}(T^\star X^{\langle n+m\rangle}).$$

Pour cela nous associons à chaque partition $\underline{e}$ de $n$ un cycle $\tau_{\underline{e}}^\star \in Z^n(T^\star X^{\langle n\rangle})$ et calculons notamment leurs produits pour $n$ et $\underline{e}$ variables. La définition de ces cycles demande certaines précautions : contrairement au cas de la caractéristique nulle, ils ne sont pas nécessairement des adhérences de conormaux de sous-variétés. En particulier, nous généralisons en caractéristique quelconque (1.4.4) une formule obtenue par Gérard Laumon en caractéristique nulle dans [Laumon 1987a, p. 330-331]. Le théorème 1.5.1, obtenu en 2020, est un cas particulier du théorème [Sawin 2021, 2.10], comme nous en a informé son auteur après notre dépôt sur arXiv.

Le calcul de ces cycles caractéristiques, joint aux propriétés différentielles du morphisme d'Abel-Jacobi $\int^n : X^{\langle n\rangle} \to \mathrm{Pic}^n_X$, permet de retrouver un théorème d'acyclicité de Pierre Deligne ([Deligne 1980]) motivé par les considérations qui vont suivre.

**0.2**. Supposons maintenant que $X$ provienne par extension des scalaires d'une courbe propre et lisse $X_1$ sur un corps fini $k_1$ de caractéristique $p > 0$, et $\mathscr{F}$ d'un faisceau étale constructible $\mathscr{F}_1 \in \mathrm{D}^{\mathrm{b}}_{\mathrm{c}}(X_1, \Lambda)$ sur $X_1$. Motivé par le programme de Langlands — et son lien avec le principe de récurrence de Deligne ([Laumon 1987b, 3.2.2.3]) — ou le formalisme de Grothendieck-Riemann-Roch ([Deligne 1987]), on souhaite exprimer, à isomorphisme *unique* près , le déterminant du complexe $\mathrm{R}\Gamma(X, \mathscr{F})$, muni d'une action du groupe de Galois de $k_1$, en termes des singularités de $\mathscr{F}$. Plus précisément, G. Laumon démontre à la fin des années 1980 la « formule du produit » ([Laumon 1987b, 3.2.1.1]) qui exprime ce $\Lambda$-module de rang un comme un produit tensoriel, indexé par les points de $X$, de $\Lambda$-modules de rang un « locaux », sous une forme qui « relève » la formule de Grothendieck-Ogg-Šafarevič ([SGA 5 X]) : la multiplication par un élément $t \in \Lambda^\times$ agissant sur $\mathscr{F}$ induit la multiplication par $t^{\chi(X,\mathscr{F})}$ sur le déterminant, une factorisation de celui-ci entraîne notamment une décomposition de la caractéristique d'Euler-Poincaré. Sa méthode repose sur la transformée de Fourier locale, introduite à cette occasion, et donc en particulier sur une réduction au cas de la droite projective. Précédemment P. Deligne avait introduit une méthode *géométrique* lui ayant permis d'établir la formule du produit, d'abord dans le cas abélien (faisceau de rang 1 [Deligne 1974]) puis dans le cas modéré ([Deligne 1980]). Nous terminons ce travail par une présentation rapide d'*une partie* de ces derniers résultats : on explique notamment pourquoi les déterminants de la cohomologie de $X$ à valeurs dans deux faisceaux étale-localement isomorphes sont canoniquement isomorphes, sous l'hypothèse qu'ils sont modérément ramifiés. (En particulier, les déterminants du Frobenius sont égaux.)

Le point de départ de la méthode de P. Deligne est l'expression, par la formule de Künneth symétrique ([SGA 4 XVII, 5.5.21]) et la formule du décalage de Quillen-Illusie ([Quillen 1968, prop. 7.21], [Illusie 1971-1972, I.4.3.2.1]), du déterminant de la cohomologie de $X$ (le facteur epsilon) comme celui de la cohomologie d'une *puissance symétrique* de $X$. Poussant par le morphisme d'Abel-Jacobi vers la jacobienne, et utilisant le fait que la structure de la cohomologie d'une variété abélienne induit (en genre $\geqslant 2$) une trivialisation canonique du déterminant à valeurs dans un faisceau lisse ([Deligne 1974, p. 70]), il est aisé (voir 3.3) d'exprimer le facteur epsilon comme le déterminant de la fibre d'un complexe de cycles évanescents : l'ingrédient manquant est le résultat d'acyclicité susmentionné pour le morphisme d'Abel-Jacobi (relativement à la puissance tensorielle symétrique externe du faisceau de départ). Ce complexe $\Phi := \Phi_{t\circ\int^n}(\mathscr{F}^{\langle n\rangle})$ dépend de deux paramètres : l'entier $n$, égal à la caractéristique d'Euler-Poincaré de $\mathscr{F}[1]$, et un germe (*arbitraire*) de fonction lisse $t : \mathrm{Pic}^n_X \dashrightarrow \mathbb{A}^1_k$ ; le point en lequel $\Phi$ est supporté dépend à la fois de $\mathscr{F}$ (uniquement à travers $S$, $r$ et les $a_s$ comme en 0.1) et, infinitésimalement, de $t$ à travers la forme différentielle $\omega_t \in H^0(X, \Omega^1_{X/k})$ naturellement associée.

## 1. Cycle caractéristique d'une puissance symétrique

### 1.1. Énoncés.

**1.1.1.** Pour toute une famille à support fini d'entiers naturels $\underline{e} = (e_i)_{i>0}$, on note $|\underline{e}|$ (resp. $\|\underline{e}\|$ la somme $\sum_i e_i$ (resp. $\sum_i ie_i$). Rappelons (A.1.2) que pour tout $n \in \mathbf{N}$ et toute courbe quasi-projective lisse $X$ sur un corps algébriquement clos $k$, on note $X^{\langle n\rangle}$ le produit symétrique $n$-ième de $X$. Pour $\underline{e}$ telle que $\|\underline{e}\| = n$, on note $\pi_{\underline{e}}$ le morphisme (fini) canonique

$$X^{\langle \underline{e}\rangle} := \prod_i X^{\langle e_i\rangle} \to X^{\langle n\rangle}$$

défini par la formule « $(D_1, \dots, D_n) \longmapsto D_1 + 2D_2 + \cdots + 3D_3 + \cdots + nD_n$ » au niveau des $k$-points. (En effet, $D_i \in X^{\langle e_i\rangle}(k)$ correspond à un cycle effectif de degré $e_i$ sur la courbe $X$.) L'image inverse $\pi^\star_{\underline{e}} TX^{\langle n\rangle}$ du fibré tangent de $X^{\langle n\rangle}$ s'identifie à l'image directe par la projection $p_2 : X \times X^{\langle \underline{e}\rangle} \to X^{\langle \underline{e}\rangle}$ de $\mathscr{O}(\mathscr{D}_1 + 2\mathscr{D}_2 + \cdots + n\mathscr{D}_n)/\mathscr{O}$ où $\mathscr{D}_i$ est le diviseur de Cartier effectif relatif de degré $e_i$ associé au facteur $X^{\langle e_i\rangle}$.

**1.1.2. Définition.** *L'inclusion évidente*

$$p_{2,\star}\left(\mathscr{O}(\mathscr{D}_1 + \mathscr{D}_2 + \cdots + \mathscr{D}_n)/\mathscr{O}\right) \subset p_{2,\star}\left(\mathscr{O}(\mathscr{D}_1 + 2\mathscr{D}_2 + \cdots + n\mathscr{D}_n)/\mathscr{O}\right)$$

*induit par passage à l'orthogonal un sous-fibré vectoriel* $\tilde{T}^\star_{\underline{e}} X^{\langle n\rangle}$ *de* $\pi^\star_{\underline{e}} T^\star X^{\langle n\rangle}$ *de rang* $n - \sum_i e_i$. *On note* $\tau^\star_{\underline{e}} X^{\langle n\rangle}$ *(ou simplement* $\tau^\star_{\underline{e}}$*) le cycle algébrique* $\pi_{\underline{e},\star}[\tilde{T}^\star_{\underline{e}} X^{\langle n\rangle}]$ *de dimension* $n$ *sur* $T^\star X^{\langle n\rangle}$ *(où, par abus de notation, on a encore noté* $\pi_{\underline{e}}$ *le morphisme canonique* $\pi^\star_{\underline{e}} T^\star X^{\langle n\rangle} \to T^\star X^{\langle n\rangle}$*).*

**1.1.3. Définition.** *En conservant les notations de la définition 1.1.2, on note* $\widehat{X^{\langle \underline{e}\rangle}}$ *le sous-schéma ouvert de* $X^{\langle \underline{e}\rangle}$ *dont les* $k$*-points correspondent aux familles de diviseurs effectifs* $(D_1, \dots, D_n)$ *sans multiplicités et à supports disjoints. (Alternativement, si on identifie* $X^{\langle \underline{e}\rangle}$ *au quotient* $X^{|\underline{e}|}/H$ *où* $H := \mathfrak{S}_{e_1} \times \cdots \times \mathfrak{S}_{e_n}$, $\widehat{X^{\langle \underline{e}\rangle}}$ *s'identifie au quotient par* $H$ *de l'ouvert invariant de* $X^{|\underline{e}|}$ *obtenu en enlevant toutes les diagonales partielles.) Notons* $X^{\langle n\rangle}_{\underline{e}}(k)$ *le sous-ensemble de* $X^{\langle n\rangle}(k)$ *formé des* $k$*-points appartenant à l'image du morphisme composé* $\widehat{X^{\langle \underline{e}\rangle}} \to X^{\langle \underline{e}\rangle} \xrightarrow{\pi_{\underline{e}}} X^{\langle n\rangle}$.

Tout $k$-point de $X^{\langle n\rangle}$ s'identifie à un diviseur effectif de degré $n$ que l'on peut écrire $D=\sum_{j=1}^{l}\mu_j x_j$ avec les $x_j$ distincts. Si pour tout $i$, on note $e_i$ le nombre d'indices $j$ tels que $\mu_j=i$, alors $\underline{e}=(e_1,\dots,e_n)$ est l'unique famille d'entiers naturels telle que $D\in X_{\underline{e}}^{\langle n\rangle}(k)$. Nous allons justifier que l'on obtient bien ainsi une stratification de $X^{\langle n\rangle}$ par des localement fermés.

Si $\mathscr{A}$ est une partition de l'ensemble $\{1,\dots,n\}$, on peut écrire $\mathscr{A}=\{A_1,\dots,A_l\}$, poser $\mu_j := \#A_j$, de sorte que $\mu_1+\dots+\mu_l=n$, puis noter pour tout $i>0$, $e_i := \#\{j,\mu_j=i\}$. On peut alors introduire la diagonale partielle $X^{\mathscr{A}}\subset X^n$ définie de façon à ce que pour tout $k$-schéma $S$, un $n$-uplet $(x_1,\dots,x_n)\in X(S)^n$ appartienne à $X^{\mathscr{A}}(S)$ si et seulement si pour tout $A\in\mathscr{A}$, si $\{j,j'\}\subset A$, alors $x_j=x_{j'}$.

Un $k$-point $D\in X^{\langle n\rangle}(k)$ appartient à $D\in X_{\underline{e}}^{\langle n\rangle}(k)$ si et seulement s'il peut s'écrire $D=\sum_j \mu_j x_j$ avec les $x_j\in X(k)$ distincts. Si on relâche la condition que les $x_j$ soient distincts, on obtient un sous-schéma fermé de $X^{\langle n\rangle}$ que l'on va noter $\overline{X_{\underline{e}}^{\langle n\rangle}}$ et qui s'identifie à l'image (schématique) du morphisme fini $\pi_{\underline{e}}\colon X^{\langle\underline{e}\rangle}\to X^{\langle n\rangle}$. On montre alors facilement qu'il existe bien un sous-schéma ouvert (dense) $X_{\underline{e}}^{\langle n\rangle}(k)$ dans le schéma $\overline{X_{\underline{e}}^{\langle n\rangle}}$ et que le complémentaire de $X_{\underline{e}}^{\langle n\rangle}(k)$ dans $\overline{X_{\underline{e}}^{\langle n\rangle}}$ est la réunion finie des fermés $\overline{X_{\underline{e}'}^{\langle n\rangle}}$ pour $\underline{e}'$ parcourant les familles d'entiers correspondant aux familles $(\mu'_1,\dots,\mu'_{l-1})$ obtenues à partir de $(\mu_1,\dots,\mu_l)$ en enlevant deux termes $\mu_j$ et $\mu_{j'}$ et en les remplaçant par un seul terme $\mu_j+\mu_{j'}$. Alternativement, le complémentaire de $X_{\underline{e}}^{\langle n\rangle}(k)$ dans $\overline{X_{\underline{e}}^{\langle n\rangle}}$ s'identifie à la réunion disjointe (ensembliste) des localement fermés $\overline{X_{\underline{e}'}^{\langle n\rangle}}$ où $\underline{e}'$ parcourt les familles d'entiers qui correspondent *via* la description précédente à des partitions $\mathscr{A}'\neq\mathscr{A}$ de $\{1,\dots,n\}$ ayant la propriété que toute partie $A'\in\mathscr{A}'$ est réunion de parties appartenant à $\mathscr{A}$. (Autrement dit, la relation d'équivalence sur l'ensemble $\{1,\dots,n\}$ associée à $\mathscr{A}'$ est strictement plus grossière que celle associée à $\mathscr{A}$.)

On dispose ainsi d'une stratification de $X^{\langle n\rangle}$ par des localement fermés $X_{\underline{e}}^{\langle n\rangle}$ pour $\underline{e}$ parcourant les familles d'entiers naturels $(e_1,\dots,e_n)$ telles que $\sum_i ie_i=n$. En outre, si on note encore $\pi_{\underline{e}}\colon X^{\langle\underline{e}\rangle}\to \overline{X_{\underline{e}}^{\langle n\rangle}}$ le morphisme fini induit, alors $\widehat{X^{\langle\underline{e}\rangle}}=\pi_{\underline{e}}^{-1}(X_{\underline{e}}^{\langle n\rangle})$.

**1.1.4. Proposition.** *Avec les notations de la définition 1.1.2, le morphisme $\widehat{X^{\langle\underline{e}\rangle}}\to X_{\underline{e}}^{\langle n\rangle}$ est un isomorphisme de schémas si les entiers $i>0$ tels que $e_i\neq 0$ sont inversibles dans le corps $k$ (et la réciproque est vraie si $X$ est non vide).*

Si $(D_1,\dots,D_n)$ est une famille de diviseurs effectifs de degrés respectifs $(d_1,\dots,d_n)$ sur $X$, alors la différentielle de $\widehat{X^{\langle\underline{e}\rangle}}\to X^{\langle n\rangle}$ en ce point s'identifie à l'application linéaire

$$\bigoplus_{i=1}^{n}\varGamma(X,\mathscr{O}(D_i)/\mathscr{O})\to\varGamma(X,\mathscr{O}(\textstyle\sum_i iD_i)/\mathscr{O})$$

obtenue en faisant la somme des applications composées suivantes, où la première flèche est la multiplication par $i$ :

$$\varGamma(X,\mathscr{O}(D_i)/\mathscr{O})\xrightarrow{i}\varGamma(X,\mathscr{O}(D_i)/\mathscr{O})\hookrightarrow\varGamma(X,\mathscr{O}(iD_i)/\mathscr{O})\hookrightarrow\varGamma(X,\mathscr{O}(\textstyle\sum_i iD_i)/\mathscr{O}).$$

Si on fait l'hypothèse que les diviseurs $D_i$ sont à supports disjoints, il est évident que cette différentielle est injective si et seulement si les entiers $i$ tels que $e_i\neq 0$ sont inversibles dans $k$ et qu'alors l'image de la différentielle s'identifie à $\varGamma(X,\mathscr{O}(\sum_i D_i)/\mathscr{O})$. On peut conclure en utilisant le lemme suivant au morphisme fini $\widehat{X^{\langle\underline{e}\rangle}}\to X_{\underline{e}}^{\langle n\rangle}$ :

**1.1.5. Lemme.** *Soit $k$ un corps algébriquement clos. Soit $f : Y \to X$ un morphisme fini entre $k$-schémas de type fini. On suppose que $f$ induit une injection $Y(k) \to X(k)$ et que pour tout $y \in Y(k)$, la différentielle $T_y Y \to T_{f(y)} X$ est injective. Alors $f$ est une immersion fermée.*

Procéder comme dans [Hartshorne 1977, Proposition 7.3 & Lemma 7.4, Chapter II].

**1.1.6. Corollaire.** *Avec les notations de la définition 1.1.2, si les entiers $i$ tels que $e_i \neq 0$ sont inversibles dans $k$, alors $X_{\underline{e}}^{\langle n\rangle}$ est lisse et $\tau_{\underline{e}}^{\star} X^{\langle n\rangle}$ est le cycle algébrique de l'adhérence du conormal de $X_{\underline{e}}^{\langle n\rangle}$ dans $X^{\langle n\rangle}$.*

La formule suivante pour les faisceaux localement constants fait partie du formulaire du §1.4, voir la proposition 1.4.8 :

**1.1.7. Théorème.** *Soit $X$ une courbe quasi-projective lisse sur un corps algébriquement clos $k$. Soit $\Lambda$ un corps fini de caractéristique inversible dans $k$. Soit $\mathscr{F}$ un faisceau localement constant de $\Lambda$-modules sur $X$. On suppose que $\mathscr{F}$ est de rang constant $r$. Pour tout entier naturel $n$, le cycle caractéristique de la puissance tensorielle symétrique externe $\mathscr{F}^{\langle n\rangle}$ (cf. A.1.3) sur $X^{\langle n\rangle}$ est donné par la formule suivante :*

$$CC(X^{\langle n\rangle}, \mathscr{F}^{\langle n\rangle}) = (-1)^n \sum_{\underline{e}, \|\underline{e}\|=n} \left( \prod_{i=1}^{n} \binom{r}{i}^{e_i} \right) \cdot \tau_{\underline{e}}^{\star}$$

*et le support singulier $SS(X^{\langle n\rangle}, \mathscr{F}^{\langle n\rangle})$ est exactement le support du cycle caractéristique.*

Pour étudier le cas de faisceaux non localement constants, il est nécessaire d'introduire d'autres cycles algébriques :

**1.1.8. Définition.** *Soit $\Delta$ un $0$-cycle effectif sur une courbe quasi-projective lisse $X$ sur un corps algébriquement clos $k$. Soit $m \in \mathbf{N}$. Soit $\underline{e} = (e_1, \dots, e_m)$ des entiers tels que $\|\underline{e}\| = m$. Notons $n := |\Delta| + m$ et $i : X^{\langle m\rangle} \to X^{\langle n\rangle}$ l'immersion fermée donnée par la formule $D \longmapsto \Delta + D$. On dispose alors d'un diagramme canonique où $p$ est lisse de dimension relative $|\Delta|$ et où $\tilde{\imath}$ est une immersion fermée :*

$$\begin{array}{ccc} & i^{\star}T^{\star}X^{\langle n\rangle} & \\ {\scriptstyle p} \swarrow & & \searrow {\scriptstyle \tilde{\imath}} \\ T^{\star}X^{\langle m\rangle} & & T^{\star}X^{\langle n\rangle} \end{array}$$

*On définit un cycle $\tau_{\Delta,\underline{e}}^{\star} X^{\langle n\rangle}$ (ou $\tau_{\Delta,\underline{e}}^{\star}$) sur $T^{\star}X^{\langle n\rangle}$ par la formule $\tau_{\Delta,\underline{e}}^{\star} X^{\langle n\rangle} := \tilde{\imath}_{\star} p^{\star} \tau_{\underline{e}}^{\star} X^{\langle m\rangle}$.*

**1.1.9. Remarque.** Si $\Delta$ est un 0-cycle effectif de degré $n$ sur une courbe quasi-projective lisse $X$ sur un corps algébriquement clos, on notera simplement $\tau_{\Delta}^{\star} := \tau_{\Delta,\underline{e}}^{\star}$ où $\underline{e}$ est la suite vide tautologique (telle que $\|\underline{e}\| = 0$). On vérifie facilement qu'alors $\tau_{\Delta}^{\star}$ est le conormal de l'immersion fermée $\mathrm{Spec}(k) \to X^{\langle n\rangle}$ donnée par $\Delta \in X^{\langle n\rangle}(k)$.

## 1.2. Image directe par les morphismes $X^{\langle n\rangle} \times X^{\langle n'\rangle} \to X^{\langle n+n'\rangle}$.

**1.2.1. Définition.** *Soit $X$ une courbe quasi-projective lisse sur un corps algébriquement clos $k$. Soit $n \in \mathbf{N}$. Un fermé de $T^{\star}X^{\langle n\rangle}$* est basique *s'il est contenu dans une réunion finie de supports de cycles de la forme $\tau_{\Delta,\underline{e}}^{\star}$ pour $|\Delta| + \|\underline{e}\| = n$. Un cycle algébrique de codimension $n$ sur $T^{\star}X^{\langle n\rangle}$* est basique *si son support l'est.*

**1.2.2. Définition.** *Soit $X$ une courbe quasi-projective lisse sur un corps algébriquement clos $k$. Soit $(n, n') \in \mathbf{N}^2$. Soit $S$ (resp. $S'$) un fermé de $T^{\star}X^{\langle n\rangle}$ (resp. $T^{\star}X^{\langle n'\rangle}$). Notons $p := \vee_{n,n'} : X^{\langle n\rangle} \times X^{\langle n'\rangle} \to$*

*$X^{\langle n+n'\rangle}$ le morphisme $(D, D') \longmapsto D+D'$ correspondant à la somme des diviseurs effectifs. On introduit le diagramme suivant :*

$$\begin{array}{ccc} & p^\star T^\star X^{\langle n+n'\rangle} & \\ {\scriptstyle dp^\vee}\swarrow & & \searrow{\scriptstyle \tilde{p}} \\ T^\star X^{\langle n\rangle} \times T^\star X^{\langle n'\rangle} & & T^\star X^{\langle n+n'\rangle} \end{array}$$

*où $dp^\vee$ est le morphisme dual de la différentielle de $p$, vue comme morphisme entre fibrés vectoriels sur $X^{\langle n\rangle} \times X^{\langle n'\rangle}$, tandis que $\tilde{p}$ est le morphisme déduit de $p$ par changement de base.*

*On note $S \,\overline{\vee}\, S' := \tilde{p}_\star(dp^\vee)^{-1}(S \times S')$.*

**1.2.3. Définition.** *Soit $X$ une courbe quasi-projective lisse sur un corps algébriquement clos $k$. Soit $(n, n') \in \mathbf{N}^2$. Soit $z \in Z^n(T^\star X^{\langle n\rangle})$. Soit $z' \in Z^{n'}(T^\star X^{\langle n'\rangle})$. Avec les notations de la définition précédente, si le pull-back $(dp^\vee)^\star(z \boxtimes z')$ est bien défini en théorie des intersections, on note :*

$$z \cdot z' := \tilde{p}_\star(dp^\vee)^\star(z \boxtimes z') \in Z^{n+n'}(T^\star X^{\langle n+n'\rangle}).$$

*(Plus précisément, si on note $G \subset P := (p^\star T^\star X^{\langle n+n'\rangle}) \times_{X^{\langle n\rangle}\times X^{\langle n'\rangle}} \big(T^\star X^{\langle n\rangle}) \times T^\star X^{\langle n\rangle}\big)$ le graphe du morphisme $dp^\vee$ et que l'on note $1 \boxtimes (z \boxtimes z')$ le cycle sur $P$ obtenu en tirant en arrière $z \boxtimes z'$, alors le produit $z \cdot z'$ est défini si et seulement si l'intersection $G \cdot (1 \boxtimes (z \boxtimes z'))$ est propre, autrement dit si les composantes irréductibles dans l'intersection ensembliste associée sont toutes de dimension $n + n'$.)*

Nous allons montrer que $\overline{\vee}$ préserve les fermés basiques, et que la structure multiplicative partiellement définie sur le groupe abélien gradué $(Z^n(T^\star X^{\langle n\rangle}))_{n\in\mathbf{N}}$ est bien définie sur les cycles basiques.

L'intérêt de définir $\overline{\vee}$ et $\cdot$ vient de la proposition suivante, qui découle de cas particuliers de résultats de Takeshi Saito sur le support singulier et le cycle caractéristique :

**1.2.4. Proposition.** *Soit $X$ une courbe quasi-projective lisse sur un corps algébriquement clos $k$. Soit $\Lambda$ un corps fini de caractéristique inversible dans $k$. Soit $(n, n') \in \mathbf{N}^2$. Soient $\mathscr{F} \in \mathrm{D}^{\mathrm{b}}_{\mathrm{c}}(X^{\langle n\rangle}, \Lambda)$, $\mathscr{F}' \in \mathrm{D}^{\mathrm{b}}_{\mathrm{c}}(X^{\langle n'\rangle}, \Lambda)$ et $\mathscr{F} \vee \mathscr{F}' := \vee_{n,n',\star}(\mathscr{F} \boxtimes \mathscr{F}') \in \mathrm{D}^{\mathrm{b}}_{\mathrm{c}}(X^{\langle n+n'\rangle}, \Lambda)$ (A.1.3). Alors,*

$$SS(\mathscr{F} \vee \mathscr{F}') \subset SS(\mathscr{F}) \,\overline{\vee}\, SS(\mathscr{F}').$$

*Si le produit de $CC(\mathscr{F})$ et de $CC(\mathscr{F}')$ est bien défini et que $X$ est projective, alors*

$$CC(\mathscr{F} \vee \mathscr{F}') = CC(\mathscr{F}) \cdot CC(\mathscr{F}').$$

En effet, on connaît le support singulier et le cycle caractéristique de $\mathscr{F} \boxtimes \mathscr{F}'$ par la formule sur les produits externes [Saito 2017a, Theorem 2.2], et on peut utiliser [Saito 2021, Theorem 2.2.5] qui permet d'étudier l'image directe par $\vee_{n,n'} : X^{\langle n\rangle} \times X^{\langle n'\rangle} \to X^{\langle n+n'\rangle}$, ce qui dans le cas du cycle caractéristique nécessite de supposer que $X$ est projective. (L'hypothèse de projectivité ne nous gênera pas vraiment dans la suite puisque l'on pourra souvent appliquer cette formule à un faisceau $\overline{\mathscr{F}}$ sur la compactification canonique $\overline{X}$ de $X$ si $X$ est seulement quasi-projective lisse.)

**1.2.5. Remarque.** Dans la proposition 1.2.4, si $\mathscr{F}'$ est le faisceau gratte-ciel de fibre $\Lambda$ supporté par un point $\Delta' \in X^{\langle n'\rangle}$, alors $\mathscr{F} \vee \mathscr{F}'$ s'identifie à l'image directe de $\mathscr{F}$ par le morphisme $X^{\langle n\rangle} \to X^{\langle n+n'\rangle}$ d'addition avec $\Delta'$, et ce morphisme est une immersion fermée. Dans ce cas particulier, la formule $CC(\mathscr{F} \vee \mathscr{F}') = CC(\mathscr{F}) \cdot CC(\mathscr{F}')$ ne fait qu'énoncer une compatibilité triviale du cycle caractéristique par rapport aux immersions fermées entre schémas lisses, et cette formule vaut bien sûr même si on suppose seulement que $X$ est quasi-projective lisse.

Pour étudier la structure multiplicative définie plus haut, il sera pratique d'utiliser une définition légèrement différente des cycles $\tau^\star_{\underline{e}}$ :

**1.2.6. Définition.** *Soit $X$ une courbe quasi-projective lisse sur un corps algébriquement clos $k$. Soit $n \in \mathbf{N}$. Soit $\mathscr{A}$ une partition d'un ensemble fini $I$ de cardinal $n$. On note $X^{\mathscr{A}}$ le sous-schéma fermé du produit $X^I$ défini par les conditions $x_i = x_j$ pour tout couple $(i,j) \in I^2$ tel qu'il existe $A \in \mathscr{A}$ tel que $\{i,j\} \subset A$. On note $\pi_{\mathscr{A}} \colon X^{\mathscr{A}} \to X^{\langle n\rangle}$ le morphisme composé évident $X^{\mathscr{A}} \to X^I \simeq X^n \to X^{\langle n\rangle}$. Considérons la courbe relative donnée par le morphisme $p_2 \colon X \times X^{\mathscr{A}} \to X^{\mathscr{A}}$. Pour tout $A \in \mathscr{A}$, on note $x_A$ la section évidente du morphisme $p_2$. Pour tout $i \in A$, on note aussi $x_i := x_A$. On peut identifier $\pi_{\mathscr{A}}^\star T X^{\langle n\rangle}$ à $p_{2,\star}\mathscr{O}(\sum_{i\in I} x_i)/\mathscr{O}$. On note $\tilde{T}^\star_{\mathscr{A}} X^{\langle n\rangle}$ le sous-fibré vectoriel de $\pi_{\mathscr{A}}^\star T^\star X^{\langle n\rangle}$ de rang $n - \#\mathscr{A}$ défini par passage à l'orthogonal de $p_{2,\star}\left(\mathscr{O}(\sum_{A\in\mathscr{A}} x_A)/\mathscr{O}\right) \subset p_{2,\star}\left(\mathscr{O}(\sum_{i\in I} x_i)/\mathscr{O}\right)$. Notons $\tau^\star_{\mathscr{A}} := \pi_{\mathscr{A},\star}[\tilde{T}^\star_{\mathscr{A}} X^{\langle n\rangle}] \in Z^n(T^\star X^{\langle n\rangle})$. Si on note $\underline{e} = (e_1, \dots, e_n)$ le $n$-uplet défini par $e_i := \#\{A \in \mathscr{A}, \#A = i\}$, il est évident que $\tau^\star_{\mathscr{A}} = \underline{e}! \cdot \tau^\star_{\underline{e}}$ où $\underline{e}! := \prod_{i=1}^n e_i!$.*

**1.2.7. Proposition.** *Soit $X$ une courbe quasi-projective lisse sur un corps algébriquement clos $k$. Supposons que $\mathscr{A}$ et $\mathscr{A}'$ sont des partitions de deux ensembles disjoints $I$ et $I'$ de cardinaux respectifs $n$ et $n'$. Notons $J(\mathscr{A},\mathscr{A}')$ l'ensemble des partitions $\mathscr{C}$ de $I \sqcup I'$ telles que pour tout $C \in \mathscr{C}$, on ait $C \cap I \in \mathscr{A} \cup \{\varnothing\}$ et $C \cap I' \in \mathscr{A}' \cup \{\varnothing\}$. Alors,*

$$\left|\tau^\star_{\mathscr{A}}\right| \,\overline{\vee}\, \left|\tau^\star_{\mathscr{A}'}\right| = \bigcup_{\mathscr{C}\in J(\mathscr{A},\mathscr{A}')} \left|\tau^\star_{\mathscr{C}}\right|.$$

*(On a noté ici $|z|$ le support d'un cycle algébrique $z$.)*

*Démonstration.* Notons $q \colon X^{\mathscr{A}} \times X^{\mathscr{A}'} \to X^{\langle n\rangle} \times X^{\langle n'\rangle}$ le morphisme évident et $\pi := p \circ q$ où $p := \vee_{n,n'} \colon X^{\langle n\rangle} \times X^{\langle n'\rangle} \to X^{\langle n+n'\rangle}$ est le morphisme canonique. On peut former le diagramme suivant.

$$\begin{array}{ccccccc}
\tilde{T}^\star_{\mathscr{A}} X^{\langle n\rangle} \times \tilde{T}^\star_{\mathscr{A}'} X^{\langle n'\rangle} & \hookrightarrow & q^\star\left(T^\star X^{\langle n\rangle} \times T^\star X^{\langle n'\rangle}\right) & \xleftarrow{q^\star(dp^\vee)} & \pi^\star T^\star X^{\langle n+n'\rangle} & \xrightarrow{\tilde{\pi}} & T^\star X^{\langle n+n'\rangle} \\
 & \searrow & & \searrow & \downarrow & & \downarrow \\
 & & & & X^{\mathscr{A}} \times X^{\mathscr{A}'} & \xrightarrow{\pi} & X^{\langle n+n'\rangle}
\end{array}$$

Par construction, $\left|\tau^\star_{\mathscr{A}}\right| \,\overline{\vee}\, \left|\tau^\star_{\mathscr{A}'}\right|$ est l'image par $\tilde{\pi}$ de l'image inverse par $q^\star(dp^\vee)$ de $\tilde{T}^\star_{\mathscr{A}} X^{\langle n\rangle} \times \tilde{T}^\star_{\mathscr{A}'} X^{\langle n'\rangle}$. Nous allons raisonner sur les $k$-points et étudier la contribution à $\left|\tau^\star_{\mathscr{A}}\right| \,\overline{\vee}\, \left|\tau^\star_{\mathscr{A}'}\right|$ de la fibre au-dessus d'un point $((x_i)_{i\in I}, (x'_{i'})_{i'\in I'}) \in (X^{\mathscr{A}} \times X^{\mathscr{A}'})(k)$, où les $x_i$ et $x'_{i'}$ sont des $k$-points de $X$, qui ont ici la propriété que si $A \in \mathscr{A}$, alors les éléments $(x_i)_{i\in A}$ sont égaux à un élément que l'on notera $x_A$ et que si $A' \in \mathscr{A}'$, alors les éléments $(x'_{i'})_{i'\in A'}$ sont égaux à un élément que l'on notera $x'_{A'}$.

Notons $D := \sum_{i\in I} x_i$, $D' := \sum_{i'\in I'} x'_{i'}$ et $E = D + D'$. Notons $T := T_D X^{\langle n\rangle} \simeq \Gamma(X, \mathscr{O}(D)/\mathscr{O})$ et $T' := T_{D'} X^{\langle n'\rangle} \simeq \Gamma(X, \mathscr{O}(D')/\mathscr{O})$. L'espace vectoriel dual $(T \times T')^\vee$ s'identifie à l'ensemble des $k$-points de la fibre de $q^\star\left(T^\star X^{\langle n\rangle} \times T^\star X^{\langle n'\rangle}\right)$ au-dessus de $((x_i)_{i\in I}, (x'_{i'})_{i'\in I'})$. L'ensemble des $k$-points de la fibre de $\tilde{T}^\star_{\mathscr{A}} X^{\langle n\rangle} \times \tilde{T}^\star_{\mathscr{A}'} X^{\langle n'\rangle}$ s'identifie à $(V \times V')^\perp$ où $V := \Gamma(X, \mathscr{O}(\sum_{A\in\mathscr{A}} x_A)/\mathscr{O}) \subset T$ et $V' := \Gamma(X, \mathscr{O}(\sum_{A'\in\mathscr{A}'} x'_{A'})/\mathscr{O}) \subset T'$.

Au point considéré, la différentielle de $p \colon X^{\langle n\rangle} \times X^{\langle n'\rangle} \times X^{\langle n+n'\rangle}$ s'identifie à l'application évidente $T \times T' \xrightarrow{d} \Gamma(X, \mathscr{O}(E)/\mathscr{O})$. Il s'agit de déterminer l'image inverse de $(V \times V')^\perp$ par $d^\vee$. Par orthogonalité, on a bien sûr $d^{\vee,-1}((V \times V')^\perp) = (d(V \times V'))^\perp$. L'image $d(V \times V')$ est facile à calculer :

$$d(V \times V') = \Gamma(X, \mathscr{O}(F)/\mathscr{O}) \quad \text{où } F := \sup\left(\sum_{A\in\mathscr{A}} x_A, \sum_{A'\in\mathscr{A}'} x'_{A'}\right).$$

En général, $F$ n'est pas nécessairement la somme $\sum_{A\in\mathscr{A}} x_A + \sum_{A'\in\mathscr{A}'} x'_{A'}$, puisque qu'il peut exister $A \in \mathscr{A}$ et $A' \in \mathscr{A}'$ tels que $x_A = x'_{A'}$. Pour tout $j \in I \sqcup I'$, on note $y_j := x_j$ si $j \in I$ et $y_j := x'_j$ si $j \in I'$.

Nous allons construire une partition $\mathscr{C} \in J(\mathscr{A},\mathscr{A}')$ qui ait la propriété que si $C \in \mathscr{C}$, alors les points $(y_j)_{j\in C}$ sont tous égaux à un point que l'on notera $y_C$, et que $F = \sum_{C\in\mathscr{C}} y_C$. Ceci montrera que la contribution à $\left|\tau^\star_{\mathscr{A}}\right| \triangledown \left|\tau^\star_{\mathscr{A}'}\right|$ de la fibre au-dessus de $((x_i)_{i\in I},(x'_{i'})_{i'\in I'})$ est contenue dans $\left|\tau^\star_{\mathscr{C}}\right|$.

Notons $Z := \{y_j, j \in J\}$. Pour tout $z \in Z$, notons $J_z := \{j \in I \sqcup I', y_j = z\}$. Par construction, $J_z$ est une réunion de parties appartenant à $\mathscr{A} \cup \mathscr{A}'$. Notons $\mathscr{A}_z := \{A \in \mathscr{A}, A \subset J_z\}$ et $\mathscr{A}'_z := \{A' \in \mathscr{A}', A' \subset J_z\}$. Il est évident que $\mathscr{A}_z$ (resp. $\mathscr{A}'_z$) est une partition de $J_z \cap I$ (resp. $J_z \cap I'$). On a aussi :

$$F = \sum_{z\in Z} \max(\#\mathscr{A}_z, \#\mathscr{A}'_z) \cdot z.$$

Fixons $z \in Z$. Notons $r_z := \min(\#\mathscr{A}_z, \#\mathscr{A}'_z)$. On peut choisir deux sous-ensembles $\mathscr{B}_z \subset \mathscr{A}_z$ et $\mathscr{B}'_z \subset \mathscr{A}'_z$ tels que $\#\mathscr{B}_z = \#\mathscr{B}'_z = r_z$. On peut aussi choisir une bijection $\varphi : \mathscr{B}_z \xrightarrow{\sim} \mathscr{B}'_z$. On définit alors :

$$\mathscr{C}_z := \{A \cup \varphi(A), A \in \mathscr{B}_z\} \cup \{A, A \in \mathscr{A}_z - \mathscr{B}_z\} \cup \{A', A' \in \mathscr{A}'_z - \mathscr{B}'_z\}.$$

Par construction, $\mathscr{C}_z \in J(\mathscr{A}_z,\mathscr{A}'_z)$ et $\#\mathscr{C}_z = \max(\#\mathscr{A}_z, \#\mathscr{A}'_z)$. En particulier, $\mathscr{C}_z$ est une partition de $J_z$. Comme $I \sqcup I' = \coprod_{z\in Z} J_z$, on peut donc définir une partition $\mathscr{C}$ de $I \sqcup I'$ en posant $\mathscr{C} := \cup_{z\in Z}\mathscr{C}_z$. On vérifie immédiatement que $\mathscr{C} \in J(\mathscr{A},\mathscr{A}')$, que $(y_C)_{C\in\mathscr{C}}$ est bien défini et que l'on a $F = \sum_{C\in\mathscr{C}} y_C$. Ceci montre l'inclusion

$$\left|\tau^\star_{\mathscr{A}}\right| \triangledown \left|\tau^\star_{\mathscr{A}'}\right| \subset \bigcup_{\mathscr{C}\in J(\mathscr{A},\mathscr{A}')} \left|\tau^\star_{\mathscr{C}}\right|.$$

L'inclusion réciproque est plus facile. Les détails sont laissés au lecteur. Cependant, le point-clef est que si $(y_j)_{j\in I\sqcup I'}$ est un $k$-point de $X^{\mathscr{C}}$ pour $\mathscr{C} \in J(\mathscr{A},\mathscr{A}')$, et que l'on note $x_i := y_i$ pour $i \in I$ et $x'_{i'} := y_{i'}$ pour $i' \in I$, ce qui permet de considérer les familles $(x_i)_{i\in I} \in X^{\mathscr{A}}(k)$ et $(x'_{i'})_{i'\in I'} \in X^{\mathscr{A}'}(k)$, alors avec les notations précédentes, on a :

$$G := \sum_{C\in\mathscr{C}} y_C \geq F := \sup\left(\sum_{A\in\mathscr{A}} x_A, \sum_{A'\in\mathscr{A}'} x'_{A'}\right)$$

ce qui induit des inclusions $\Gamma(X,\mathscr{O}(G)/\mathscr{O})^\perp \subset \Gamma(X,\mathscr{O}(F)/\mathscr{O})^\perp \subset T^\star_E X^{\langle n+n'\rangle}$. □

**1.2.8. Corollaire.** *Soit $X$ une courbe quasi-projective lisse sur un corps algébriquement clos $k$. Soit $(n,n') \in \mathbf{N}^2$.*

- *Si $S \subset T^\star X^{\langle n\rangle}$ et $S' \subset T^\star X^{\langle n'\rangle}$ sont des fermés basiques, alors le fermé $S \triangledown S' \subset T^\star X^{\langle n+n'\rangle}$ aussi;*
- *Si $z \in Z^n(T^\star X^{\langle n\rangle})$ et $z' \in Z^{n'}(T^\star X^{\langle n'\rangle})$ sont des cycles basiques, alors le produit $z \cdot z' \in Z^{n+n'}(T^\star X^{\langle n+n'\rangle})$ est bien défini et est un cycle basique.*

Nous déterminerons plus loin une formule pour les produits $\tau^\star_{\underline{e}} \cdot \tau^\star_{\underline{e}'}$. Pour le moment, on peut seulement énoncer quelques formules triviales, comme $\tau^\star_\Delta \cdot \tau^\star_{\Delta'} = \tau^\star_{\Delta+\Delta'}$, ou $\tau^\star_\Delta \cdot \tau^\star_{\underline{e}} = \tau^\star_{\Delta,\underline{e}}$.

**1.2.9. Définition.** *Soit $X$ une courbe quasi-projective lisse sur un corps algébriquement clos $k$. Soit $\Lambda$ un corps fini de cardinal inversible dans $k$. Soit $\mathscr{F} \in \mathsf{D}^{\mathsf{b}}_{\mathsf{c}}(X,\Lambda)$. On note :*

$$S_{\mathscr{F}} := \sum_{n\in\mathbf{N}} CC(X^{\langle n\rangle}, \mathscr{F}^{\langle n\rangle}) \in \prod_{n\in\mathbf{N}} Z^n(T^\star X^{\langle n\rangle}).$$

*On note $\overline{S}_{\mathscr{F}} := (SS(X^{\langle n\rangle}, \mathscr{F}^{\langle n\rangle}))_{n\in\mathbf{N}}$.*

Si $\overline{S} = (S_n)_{n\in\mathbf{N}}$ et $\overline{T} = (T_n)_{n\in\mathbf{N}}$ sont des familles de fermés dans les $T^\star X^{\langle n\rangle}$ pour tout $n \in \mathbf{N}$, on note :

$$\overline{S} \mathbin{\overline{\vee}} \overline{T} := \Big( \bigcup_{p+q=n} S_p \mathbin{\overline{\vee}} T_q \Big)_{n\in\mathbf{N}}.$$

(On écrira $\overline{S} \subset \overline{T}$ si pour tout $n \in \mathbf{N}$, $S_n \subset T_n$.) On étend aussi de façon évidente le produit de la définition 1.2.3 aux « séries formelles » appartenant à $\prod_{n\in\mathbf{N}} Z^n(T^\star X^{\langle n\rangle})$.

**1.2.10. Proposition.** *Soit $X$ une courbe quasi-projective lisse sur un corps algébriquement clos $k$. Soit $\Lambda$ un corps fini de cardinal inversible dans $k$. Supposons que $\mathscr{F}'$ et $\mathscr{F}''$ soient deux objets de $\mathsf{D}^{\mathsf{b}}_{\mathsf{c}}(X,\Lambda)$. Alors,*

$$\overline{S}_{\mathscr{F}'\oplus\mathscr{F}''} \subset \overline{S}_{\mathscr{F}'} \mathbin{\overline{\vee}} \overline{S}_{\mathscr{F}''}.$$

*Si $X$ est projective lisse et que le produit $S_{\mathscr{F}'} \cdot S_{\mathscr{F}}$ est défini, alors*

$$S_{\mathscr{F}'\oplus\mathscr{F}''} = S_{\mathscr{F}'} \cdot S_{\mathscr{F}''}.$$

Pour tout $n \in \mathbf{N}$, on a un isomorphisme canonique :

$$(\mathscr{F}' \oplus \mathscr{F}'')^{\langle n\rangle} \simeq \oplus_{p+q=n} \mathscr{F}'^{\langle p\rangle} \vee \mathscr{F}''^{\langle q\rangle}.$$

La proposition résulte alors immédiatement de la proposition 1.2.4.

**1.2.11. Proposition.** *Soit $X$ une courbe quasi-projective lisse sur un corps algébriquement clos $k$. Soit $\Lambda$ un corps fini de cardinal inversible dans $k$. Soit $\mathscr{F}' \to \mathscr{F} \to \mathscr{F}'' \to \mathscr{F}'[1]$ un triangle distingué dans $\mathsf{D}^{\mathsf{b}}_{\mathsf{c}}(X,\Lambda)$. On suppose que $\overline{S}_{\mathscr{F}'}$ et $\overline{S}_{\mathscr{F}''}$ sont basiques. Alors, $\overline{S}_{\mathscr{F}} \subset \overline{S}_{\mathscr{F}'} \mathbin{\overline{\vee}} \overline{S}_{\mathscr{F}''}$ est basique.*

*Si $X$ est projective, on a de plus :*

$$S_{\mathscr{F}} = S_{\mathscr{F}'} \cdot S_{\mathscr{F}''}.$$

Quitte à remplacer $\mathscr{F}$ par un complexe qui lui est quasi-isomorphe, il existe pour tout $n \in \mathbf{N}$, une filtration sur le complexe $\mathscr{F}^{\langle n\rangle}$ dont le gradué est $(\mathscr{F}' \oplus \mathscr{F}'')^{\langle n\rangle}$. Au niveau des supports singuliers, on a alors bien sûr :

$$\overline{S}_{\mathscr{F}} \subset \overline{S}_{\mathscr{F}'\oplus\mathscr{F}''} \subset \overline{S}_{\mathscr{F}'} \mathbin{\overline{\vee}} \overline{S}_{\mathscr{F}''}.$$

On obtient aussi la formule suivante au niveau des cycles caractéristiques :

$$S_{\mathscr{F}} = S_{\mathscr{F}'\oplus\mathscr{F}''}.$$

Si $X$ est projective, on a de plus $S_{\mathscr{F}'\oplus\mathscr{F}''} = S_{\mathscr{F}'} \cdot S_{\mathscr{F}''}$.

**1.2.12. Remarque.** Soit $S = \sum_{n\in\mathbf{N}} s_n$ avec $s_n \in Z^n(T^\star X^{\langle n\rangle})$ pour tout $n \in \mathbf{N}$. Si $S$ est basique et que $s_0 = 1$, alors $S$ est inversible et son inverse donné par la formule $S^{-1} = \sum_{n\in\mathbf{N}} x^n$ pour $x := 1 - S$ est basique.

**1.2.13. Proposition.** *Soit $X$ une courbe quasi-projective sur un corps algébriquement clos $k$. Soit $\Lambda$ un corps fini de cardinal inversible sur $k$. Soit $\mathscr{F} \in \mathsf{D}^{\mathsf{b}}_{\mathsf{c}}(X,\Lambda)$.*

*Alors, $\overline{S}_{\mathscr{F}}$ est basique si et seulement si $\overline{S}_{\mathscr{F}[1]}$ est basique. Si c'est le cas et que l'on suppose de plus que $X$ est projective, alors $S_{\mathscr{F}[1]} = S_{\mathscr{F}}^{-1}$.*

*Démonstration.* On dispose du triangle distingué canonique $\mathscr{F} \to 0 \to \mathscr{F}[1] \to \mathscr{F}[1]$. D'après ce qui précède, il existe un complexe $\mathscr{Z}$ quasi-isomorphe à 0 tel que pour tout $n \in \mathbf{N}$, $\mathscr{Z}^{\langle n\rangle}$ soit muni d'une filtration dont le gradué s'identifie à $\oplus_{p+q=n} \mathscr{F}^{\langle p\rangle} \vee \mathscr{F}[1]^{\langle q\rangle}$.

Supposons que pour tout $n \in \mathbf{N}$, $SS(\mathscr{F}^{\langle n\rangle})$ est basique. Montrons par récurrence sur $n \in \mathbf{N}$ que $SS(\mathscr{F}[1]^{\langle n\rangle})$ est basique. C'est évident si $n = 0$ (ou si $n = 1$) Supposons $n \geq 1$ et que le résultat est connu jusqu'à $n-1$. Comme $\mathscr{F}[1]^{\langle n\rangle}$ est un des objets « extrêmes » dans le gradué de la filtration mentionnée plus haut, il est manifeste que $\mathscr{F}[1]^{\langle n\rangle}$ appartient à la sous-catégorie triangulée de $\mathrm{D}^{\mathrm{b}}_{\mathrm{c}}(X,\Lambda)$ engendrée par les objets $\mathscr{F}^{\langle p\rangle}\vee\mathscr{F}[1]^{\langle n-p\rangle}$ pour $p\in\{1,\dots,n\}$, et ceux-là ont un support singulier basique. On en déduit que $SS(\mathscr{F}[1]^{\langle n\rangle})$ est basique. Ainsi, $\overline{S}_{\mathscr{F}}$ basique implique $\overline{S}_{\mathscr{F}[1]}$ basique, et la réciproque se montre de la même manière.

Si on suppose de plus que $X$ est projective, l'hypothèse que $\overline{S}_{\mathscr{F}}$ et $\overline{S}_{\mathscr{F}[1]}$ sont basiques permet d'appliquer le résultat de la proposition précédente au triangle distingué $\mathscr{F}\to 0\to\mathscr{F}[1]\to\mathscr{F}[1]$, d'où la formule $S_{\mathscr{F}}\cdot S_{\mathscr{F}[1]}=1$. □

**1.2.14. Remarque.** Il résulte des propositions 1.2.10, 1.2.11 et 1.2.13 que la sous-catégorie de $\mathrm{D}^{\mathrm{b}}_{\mathrm{c}}(X,\Lambda)$ formée des $\mathscr{F}$ tels que $\overline{S}_{\mathscr{F}}$ soit basique est une sous-catégorie triangulée. En outre, si on suppose que $X$ est projective, l'application $\mathscr{F}\longmapsto S_{\mathscr{F}}$ induit un morphisme de groupes du groupe de Grothendieck de la sous-catégorie triangulée considérée vers le groupe multiplicatif des séries formelles basiques de terme constant 1.

## 1.3. Détermination de la structure d'algèbre.

**1.3.1. Théorème.** *Soit $X$ une courbe quasi-projective lisse sur un corps algébriquement clos $k$. Supposons que $\mathscr{A}$ et $\mathscr{A}'$ sont des partitions de deux ensembles disjoints $I$ et $I'$ de cardinaux respectifs $n$ et $n'$. Alors,*

$$\tau^{\star}_{\mathscr{A}}\cdot\tau^{\star}_{\mathscr{A}'}=\sum_{\mathscr{C}\in J(\mathscr{A},\mathscr{A}')}\tau^{\star}_{\mathscr{C}}.$$

*(L'ensemble de partitions $J(\mathscr{A},\mathscr{A}')$ a été défini à la proposition 1.2.7.)*

On commence par le cas particulier de la droite affine, dont la démonstration nous occupera jusqu'au lemme 1.3.7 :

**1.3.2. Proposition.** *L'énoncé du théorème 1.3.1 est vrai dans le cas particulier $X=\mathbf{A}^1_k$.*

On notera $X:=\mathbf{A}^1_k$ pendant toute la démonstration de cette proposition. Il s'agit d'étudier la différentielle de $p:=\vee_{n,n'}: X^{\langle n\rangle}\times X^{\langle n'\rangle}\to X^{\langle n+n'\rangle}$, qui est un morphisme entre fibrés vectoriels sur $X^{\langle n\rangle}\times X^{\langle n'\rangle}$. On va s'intéresser à l'image inverse de ce morphisme par le morphisme évident $X^{\mathscr{A}}\times X^{\mathscr{A}'}\to X^{\langle n\rangle}\times X^{\langle n'\rangle}$. Notons $B:=B_{\mathscr{A},\mathscr{A}'}$ l'anneau de $X^{\mathscr{A}}\times X^{\mathscr{A}'}$ que l'on identifie à l'anneau des polynômes en les variables $(x_A)_{A\in\mathscr{A}}$ et $(x_{A'})_{A'\in\mathscr{A}'}$. Si $i\in A\in\mathscr{A}$, on posera $x_i:=x_A\in B$ et si $i'\in A'\in\mathscr{A}'$, $x_{i'}:=x_{A'}\in B$. On introduit plusieurs polynômes dans $B[T]$ :

$$\pi_I:=\prod_{i\in I}(T-x_i)\qquad\qquad \pi_{I'}:=\prod_{i'\in I'}(T-x_{i'}).$$

Si $\pi\in B[T]$ est un polynôme unitaire (par exemple $\pi_I$, $\pi_{I'}$ ou $\pi_I\pi_{I'}$), on notera $B[T,\pi^{-1}]$ la $B$-algèbre obtenue en inversant $\pi$ (qui n'est pas un diviseur de zéro) et $\pi^{-1}B[T]$ le sous-$B$-module de $B[T,\pi^{-1}]$ formé des fractions $\frac{P}{\pi}$ avec $P\in B[T]$.

D'après [Polishchuk 2003, §16.1 et prop. 19.2], on peut identifier les sections globales du fibré vectoriel $TX^{\langle n+n'\rangle}\times_{X^{\langle n+n'\rangle}}(X^{\mathscr{A}}\times X^{\mathscr{A}'})$ au $B$-module libre $((\pi_I\pi_{I'})^{-1}B[T])/B[T]$. Plus précisément on peut identifier la différentielle $(dp)_{|X^{\mathscr{A}}\times X^{\mathscr{A}'}}$ au morphisme évident de $B$-modules :

$$(\pi_I^{-1}B[T])/B[T]\oplus(\pi_{I'}^{-1}B[T])/B[T]\to((\pi_I\pi_{I'})^{-1}B[T])/B[T].$$

Notons $\alpha\colon B[T]\to B((U))$ le morphisme de $B$-algèbre défini par $\alpha(T):=U^{-1}$. Autrement dit, $\alpha(f)=f(U^{-1})$. Si $\pi\in B[T]$ est unitaire de degré $d$, alors $\alpha(\pi)=\pi(U^{-1})=U^{-d}\cdot u$ où $u$ est un élément inversible de $B[[U]]$. En particulier, $\alpha(\pi)$ est inversible dans $B((U))$, ce qui permet d'étendre $\alpha$ en un morphisme de $B$-algèbres $\alpha\colon B[T,\pi^{-1}]\to B((U))$. Ce morphisme $\alpha$ induit un morphisme de $B$-modules $\overline{\alpha}\colon(\pi^{-1}B[T])/B[T]\to B((U))/B[U^{-1}]$. La multiplication par $U^{-1}$ induit un isomorphisme $B((U))/B[U^{-1}]\xrightarrow{\sim}B((U))/U^{-1}B[U^{-1}]\simeq B[[U]]$. On note $\beta\colon(\pi^{-1}B[T])/B[T]\xrightarrow{\alpha}B((U))/B[U^{-1}]\xrightarrow{\sim}B((U))/U^{-1}B[U^{-1}]\simeq B[[U]]$ le morphisme composé :

$$\begin{array}{ccc}(\pi^{-1}B[T])/B[T] & \xrightarrow{\beta} & B[[U]] \\ [f] & \longmapsto & U^{-1}f(U^{-1})\end{array}$$

On note que si $\lambda \in B$ (et que $\pi(\lambda) = 0$), on a $\beta(\frac{1}{T-\lambda}) = \frac{U^{-1}}{U^{-1}-\lambda} = \frac{1}{1-\lambda U} = \sum_{s=0}^{\infty} \lambda^s U^s$. Plus généralement, si $\pi = \prod_{j\in J}(T-\lambda_j) \in B[T]$, alors l'image de $(\pi^{-1}B[T])/B[T]$ par $\beta$ est l'ensemble $\mathcal{U}_\pi$ des séries formelles de la forme

$$\frac{\sum_{s=0}^{\#J-1} b_s U^s}{\prod_{j\in J}(1-\lambda_j U)} \in B[[U]]$$

où pour tout $s \in \{0,\dots,\#J-1\}$, $b_s \in B$.

En particulier, pour $d = n+n'$, $\beta$ permet d'identifier à $\mathcal{U}_{\pi_I \pi_{I'}}$ à l'image inverse sur $X^{\mathcal{A}} \times X^{\mathcal{A}'}$ du fibré tangent de $X^{\langle n+n' \rangle}$. L'image inverse sur $X^{\mathcal{A}} \times X^{\mathcal{A}'}$ de la différentielle de $\vee_{n,n'}$ s'identifie alors au morphisme évident d'addition :

$$\mathcal{U}_{\pi_I} \oplus \mathcal{U}_{\pi_{I'}} \to \mathcal{U}_{\pi_I \pi_{I'}}.$$

Pour étudier le fibré cotangent, on utilise la dualité suivante :

**1.3.3. Définition.** *Considérons l'application $B$-bilinéaire $\langle\cdot,\cdot\rangle : B[[U]] \times B[T] \to B$ définie de telle façon que si $f = \sum_i a_i U^i$ et $P = \sum_{i=0}^d b_i T^i$, alors $\langle f, P\rangle = \sum_{i=0}^d a_i b_i$.*

Si pour tout $B$-module $M$, on note $M^\vee := \mathrm{Hom}_B(M,B)$, alors cette dualité induit un isomorphisme $B[[T]] \xrightarrow{\sim} B[T]^\vee$.

**1.3.4. Lemme.** *Soit $f \in B[[U]]$. Soit $P \in B[T]$. Pour tout $\lambda \in B$, on a :*

$$\langle f, (T-\lambda)P\rangle = \langle (1-\lambda U)f, TP\rangle .$$

*Si $(\lambda_j)_{j\in J}$ est une famille finie d'éléments de $B$, on a :*

$$\left\langle \frac{f}{\prod_{j\in J}(1-\lambda_j U)}, \prod_{j\in J}(T-\lambda_j)\cdot P \right\rangle = \left\langle f, T^{\#J}P\right\rangle .$$

On dispose de la relation évidente $\langle Uf, TP\rangle = \langle f, P\rangle$. On en déduit :

$$\begin{aligned}\langle f,(T-\lambda)P\rangle &= \langle f,TP\rangle - \lambda\langle f,P\rangle \\ &= \langle f,TP\rangle - \lambda\langle Uf,TP\rangle \\ &= \langle f-\lambda Uf,TP\rangle \\ &= \langle (1-\lambda U)f,TP\rangle\end{aligned}$$

Pour tout $g \in B[[U]]$, en appliquant la relation précédente à $f := \frac{g}{1-\lambda U}$, on obtient :

$$\left\langle \frac{g}{1-\lambda U}, (T-\lambda)P \right\rangle = \langle g, TP\rangle .$$

La dernière formule énoncée dans le lemme se déduit alors immédiatement par récurrence sur $\#J$.

Il résulte du lemme que si $\pi = \prod_{j\in J}(T-\lambda_j) \in B[T]$, alors $\mathcal{U}_\pi$ est orthogonal à l'idéal $\pi B[T]$. Plus précisément, on a un isomorphisme induit $B[T]/(\pi) \xrightarrow{\sim} \mathcal{U}_\pi^\vee$. On en déduit que la duale $dp^\vee$ de la différentielle de $p := \vee_{n,n'}$ correspond, après image inverse à $X^{\mathscr{A}} \times X^{\mathscr{A}'}$ au morphisme évident de $B$-modules libres :

$$d^{\vee}_{I,I'} : B[T]/(\pi_I \pi_{I'}) \to B[T]/(\pi_I) \times B[T]/(\pi_{I'}).$$

Notons $\pi_{\mathscr{A}} := \prod_{A \in \mathscr{A}} (T - x_A) \in B[T]$ et $\pi_{\mathscr{A}'} := \prod_{A' \in \mathscr{A}'} (T - x_{A'}) \in B[T]$. D'après ce qui précède, il est évident que l'orthogonal (dans $\mathscr{U}^{\vee}_{\pi_I}$) de $\mathscr{U}_{\pi_{\mathscr{A}}} \subset \mathscr{U}_{\pi_I}$ s'identifie à $\pi_{\mathscr{A}} B[T]/\pi_I B[T] \subset B[T]/\pi_I B[T]$. Le même résultat vaut bien sûr pour l'orthogonal de $\mathscr{U}_{\mathscr{A}'}$.

On en déduit que dans l'image inverse de $T^{\star} X^{\langle n \rangle} \times T^{\star} X^{\langle n' \rangle}$ sur $X^{\mathscr{A}} \times X^{\mathscr{A}'}$, le sous-fibré vectoriel $\tilde{T}^{\star}_{\mathscr{A}} X^{\langle n \rangle} \times \tilde{T}^{\star}_{\mathscr{A}'} X^{\langle n' \rangle}$ correspond au sous-$B$-module (facteur direct)

$$(\pi_{\mathscr{A}} B[T]/\pi_I B[T]) \times (\pi_{\mathscr{A}'} B[T]/\pi_{I'} B[T]) \subset B[T]/(\pi_I) \times B[T]/(\pi_{I'}).$$

Ce sous-module est le noyau du morphisme surjectif

$$B[T]/(\pi_I) \times B[T]/(\pi_{I'}) \to B[T]/(\pi_{\mathscr{A}}) \times B[T]/(\pi_{\mathscr{A}'}).$$

Pour démontrer la proposition 1.3.2, nous allons utiliser la définition suivante :

**1.3.5. Définition.** *Soit $S$ un schéma de type fini et lisse sur un corps algébriquement clos $k$. Soit $\varphi : E \to F$ un morphisme entre fibrés vectoriels de rangs respectifs $r_E$ et $r_F$ sur $S$ ; on considère ici $E$ et $F$ comme des objets géométriques : des $X$-schémas munis d'une action du schéma en anneaux $\mathbf{A}^1_S$. Dans $E \times_S F$, on définit deux sous-schémas fermés : le graphe $G$ de $\varphi$ (défini par les équations vectorielles « $f = \varphi(e)$ ») et $E \times \{0\}$ (défini par l'équation $f = 0$). Si l'intersection de ces deux sous-schémas fermés de $E \times_S F$ est propre, on peut identifier le produit d'intersection à un cycle sur $E$ que l'on note $\mathscr{K}[\varphi] \in Z^{r_F}(E)$. On dira que $\mathscr{K}[\varphi]$ est le noyau de $\varphi$ au sens de la théorie de l'intersection.*

*Plus généralement, si $F' \subset F$ est un sous-fibré vectoriel (localement facteur direct), on peut définir de façon similaire l'image inverse de $F'$ par $\varphi$ au sens de la théorie de l'intersection, et elle coïncide avec $K[p \circ \varphi]$ où $p : F \to F/F'$.*

On notera que ces constructions sont compatibles aux changements de base plats $S' \to S$.

En appliquant cette construction au $k$-schéma $S := X^{\mathscr{A}} \times X^{\mathscr{A}'}$, il vient que pour démontrer la proposition 1.3.2, il suffit d'établir le lemme suivant :

**1.3.6. Lemme.** *Soit $k$ un corps algébriquement clos. On note $X := \mathbf{A}^1_k$. Supposons que $\mathscr{A}$ et $\mathscr{A}'$ sont des partitions de deux ensembles finis disjoints $I$ et $I'$ de cardinaux respectifs $n$ et $n'$. Notons $B := B_{\mathscr{A},\mathscr{A}'}$ l'anneau du schéma affine $X^{\mathscr{A}} \times X^{\mathscr{A}'}$. Notons $\varphi : B[T]/(\pi_I \pi_{I'}) \to B[T]/(\pi_{\mathscr{A}}) \times B[T]/(\pi_{\mathscr{A}'})$ le morphisme canonique de $B$-modules libres, que l'on identifie à un morphisme de fibrés vectoriels sur $X^{\mathscr{A}} \times X^{\mathscr{A}'}$. Alors, on a l'égalité de cycles sur $(X^{\mathscr{A}} \times X^{\mathscr{A}'}) \times_{X^{\langle n \rangle}} T^{\star} X^{\langle n+n' \rangle}$ :*

$$\mathscr{K}[\varphi] = \sum_{\mathscr{C} \in J(\mathscr{A},\mathscr{A}')} [\tilde{T}^{\star}_{\mathscr{C}} X^{\langle n+n' \rangle}].$$

D'après la proposition 1.2.7, on sait déjà que le cycle $\mathscr{K}[\varphi]$ est bien défini, puisque l'intersection considérée est propre. Remarquons que le morphisme vers le premier facteur

$$B[T]/(\pi_I \pi_{I'}) \to B[T]/(\pi_{\mathscr{A}})$$

est surjectif, de noyau $(\pi_{\mathscr{A}} B[T]/(\pi_I \pi_{I'} B[T])$. En utilisant des propriétés très élémentaires du produit d'intersection, on obtient :

$$\mathscr{K}[\varphi] = \mathscr{K}\left[(\pi_{\mathscr{A}} B[T])/(\pi_I \pi_{I'}) \to B[T]/(\pi_{\mathscr{A}'})\right].$$

Plus généralement, pour tout entier $d \geq \#\mathscr{A} + \#\mathscr{A}'$, nous allons calculer

$$\mathscr{K}\left[\pi_{\mathscr{A}} B[T]_{<d-\#\mathscr{A}} \to B[T]/(\pi_{\mathscr{A}'})\right]$$

où pour tout entier $s$, $B[T]_{<s}$ désigne le sous-module des polynômes de degré $< s$ :

**1.3.7. Lemme.** *Avec les notations précédentes, si $\mathscr{A}$ et $\mathscr{A}'$ sont des partitions de deux ensembles finis disjoints $I$ et $I'$ et que l'on note $B := B_{\mathscr{A},\mathscr{A}'}$ l'anneau du schéma affine $X^{\mathscr{A}} \times X^{\mathscr{A}'}$ (où $X := \mathbf{A}^1_k$), on a l'égalité suivante de cycles pour tout $d \geq \#\mathscr{A} + \mathscr{A}'$ :*

$$\mathscr{K}\left[\pi_{\mathscr{A}} B[T]_{<d-\#\mathscr{A}} \to B[T]/(\pi_{\mathscr{A}'})\right] = \sum_{\mathscr{C} \in J(\mathscr{A},\mathscr{A}')} \left[\pi_{\mathscr{C}} P^{d-\#\mathscr{C}}_{|X^{\mathscr{C}}}\right]$$

*où pour tout $\mathscr{C} \in J(\mathscr{A}, \mathscr{A}')$ on identifie $X^{\mathscr{C}}$ à un sous-schéma fermé de $X^{\mathscr{A}} \times X^{\mathscr{A}'}$, où pour tout entier $s$, on a noté $P^s$ le fibré vectoriel de rang $s$ sur $X^{\mathscr{A}} \times X^{\mathscr{A}'}$ correspondant au $B$-module $B[T]_{<s}$ et où de façon semblable à ce qui précède on a noté $\pi_{\mathscr{C}} := \prod_{C \in \mathscr{C}} (T - x_C) \in B_{\mathscr{C}}[T]$, avec $B_{\mathscr{C}}$ l'anneau du schéma affine $X^{\mathscr{C}}$.*

On démontre ce lemme par récurrence sur $\#\mathscr{A}'$. Le résultat est évident si $\mathscr{A}'$ (et donc $I'$) est vide. Supposons que $\mathscr{A}'$ soit non vide et choisissons $A'_0 \in \mathscr{A}'$. Notons $\mathscr{A}'' := \mathscr{A}' - \{A'_0\}$, qui est une partition de $I'' := I' - \{A'_0\}$. On observe que l'on dispose d'un isomorphisme de $B$-modules :

$$B[T]/(\pi_{\mathscr{A}'}) \xrightarrow{\sim} B[T]/(\pi_{\mathscr{A}''}) \times B$$

qui à la classe de $P = Q\pi_{\mathscr{A}''} + R$ avec $\deg R < \#\mathscr{A}''$ associe le couple $([P], Q(x_{A'_0})) = ([R], Q(x_{A'_0}))$. En utilisant cet isomorphisme, on est ramené au calcul du noyau au sens de la théorie de l'intersection du morphisme composé :

$$\varphi' : \pi_{\mathscr{A}} B[T]_{<d-\#\mathscr{A}} \to B[T]/(\pi_{\mathscr{A}'}) \times B.$$

On va calculer $\mathscr{K}[\varphi']$ en déterminant d'abord le noyau du morphisme vers le premier facteur, puis en intersectant avec l'hypersurface $\mathscr{H}$ dont l'équation correspond au morphisme vers le deuxième facteur.

En appliquant l'hypothèse de récurrence à $\mathscr{A}''$ et en utilisant la lissité du morphisme canonique $X^{\mathscr{A}} \times X^{\mathscr{A}'} \to X^{\mathscr{A}} \times X^{\mathscr{A}''}$, on obtient :

$$\mathscr{K}[\pi_{\mathscr{A}} B[T]_{<d-\#\mathscr{A}} \to B[T]/(\pi_{\mathscr{A}''})] = \sum_{\mathscr{D} \in J(\mathscr{A},\mathscr{A}'')} \left[\pi_{\mathscr{D}} P^{d-\#\mathscr{D}}_{|X^{\mathscr{D} \sqcup \{A'_0\}}}\right].$$

Il s'agit maintenant de calculer l'intersection avec l'hypersurface $\mathscr{H}$ de chacun des termes $\left[\pi_{\mathscr{D}} P^{d-\#\mathscr{D}}_{|X^{\mathscr{D} \sqcup \{A'_0\}}}\right]$ pour $\mathscr{D} \in J(\mathscr{A}, \mathscr{A}'')$.

Soit $\mathscr{D} \in J(\mathscr{A}, \mathscr{A}'')$. Notons $\tilde{\mathscr{D}} := \mathscr{D} \sqcup \{A'_0\}$ et $B_{\tilde{\mathscr{D}}} := k[(x_D)_{D \in \tilde{\mathscr{D}}}]$. Le produit d'intersection du cycle $\left[\pi_{\mathscr{D}} P^{d-\#\mathscr{D}}_{|X^{\tilde{\mathscr{D}}}}\right]$ avec $[\mathscr{H}]$ peut s'interpréter comme le noyau au sens de la théorie de l'intersection du morphisme de $B_{\tilde{\mathscr{D}}}$-modules libres :

$$\psi_{\mathscr{D}} : \pi_{\mathscr{D}} B_{\tilde{\mathscr{D}}}[T]_{<d-\#\mathscr{D}} \to B_{\tilde{\mathscr{D}}}$$

donné par la formule $\psi_{\mathscr{D}}(\pi_{\mathscr{D}} P) = \left(\frac{\pi_{\mathscr{D}}}{\pi_{\mathscr{A}''}} P\right)(x_{A'_0})$. Notons $\overline{\mathscr{D}} := \{A \in \mathscr{A}, A \in \mathscr{D}\}$. Il est immédiat que

$$\frac{\pi_{\mathscr{D}}}{\pi_{\mathscr{A}''}} = \pi_{\overline{\mathscr{D}}} := \prod_{A \in \overline{\mathscr{D}}} (T - x_A) \in B_{\tilde{\mathscr{D}}}[T].$$

Ainsi, pour tout $P \in B_{\tilde{\mathscr{D}}}[T]_{<d-\#\mathscr{D}}$, on a :

$$\psi_{\mathscr{D}}(\pi_{\mathscr{D}} P) = P(x_{A'_0}) \cdot \prod_{A \in \overline{\mathscr{D}}} (x_{A'_0} - x_A).$$

Cette factorisation implique que $\mathscr{K}[\psi_{\mathscr{D}}]$ se décompose en la somme des cycles des hypersurfaces de $\pi_{\mathscr{D}} P^{d-\#\mathscr{D}}_{|X^{\widetilde{\mathscr{D}}}}$ définies par les équations $x_{A'_0} = x_A$ pour $A \in \overline{\mathscr{D}}$ et du cycle de l'hypersurface définie par l'annulation de $P(x_{A'_0})$.

Pour tout $A \in \overline{\mathscr{D}}$, l'hypersurface définie par l'équation $x_{A'_0} = x_A$ est $\pi_{\mathscr{C}} P^{d-\#\mathscr{C}}_{|X^{\mathscr{C}}}$ où $\mathscr{C}$ est la partition obtenue à partir de $\widetilde{\mathscr{D}}$ en regroupant $A$ et $A'_0$, c'est-à-dire $\mathscr{C} = \mathscr{A}'' \sqcup \mathscr{A} - \{A\} \sqcup \{A \cup A'_0\}$. Par ailleurs, il est évident que l'hypersurface définie par l'annulation de $P(x_{A'_0})$ est $\pi_{\widetilde{\mathscr{D}}} P^{d-\#\widetilde{\mathscr{D}}}_{|X^{\widetilde{\mathscr{D}}}}$. Ainsi, les hypersurfaces obtenues sont exactement les $\pi_{\mathscr{C}} P^{d-\#\mathscr{C}}_{|X^{\mathscr{C}}}$ où $\mathscr{C}$ parcourt l'ensemble des partitions $\mathscr{C} \in J(\mathscr{A}, \mathscr{A}')$ de $I \sqcup I'$ dont la partition induite sur l'ensemble $(I \sqcup I') - A'_0$ soit $\mathscr{D}$. Par conséquent, on a l'identité suivante pour tout $\mathscr{D} \in J(\mathscr{A}, \mathscr{A}'')$ :

$$\mathscr{K}[\psi_{\mathscr{D}}] = \sum_{\substack{\mathscr{C} \in J(\mathscr{A}, \mathscr{A}') \\ \mathscr{C} \text{ induisant } \mathscr{D}}} \left[\pi_{\mathscr{C}} P^{d-\#\mathscr{C}}_{|X^{\mathscr{C}}}\right].$$

En faisant la somme sur tous les $\mathscr{D} \in J(\mathscr{A}, \mathscr{A}'')$, on obtient finalement la relation voulue :

$$\mathscr{K}[\varphi'] = \sum_{\mathscr{C} \in J(\mathscr{A}, \mathscr{A}')} \left[\pi_{\mathscr{C}} P^{d-\#\mathscr{C}}_{|X^{\mathscr{C}}}\right].$$

Ceci achève la démonstration de la proposition 1.3.2.

**1.3.8. Définition.** *Soit $X \to S$ un morphisme entre schémas lisses de type fini sur un corps algébriquement clos $k$. Supposons que $z$ et $z'$ soient deux cycles algébriques sur $X$. Soit $s \in S(k)$. On dit que $z$ et $z'$ coïncident au voisinage de $s$ s'il existe un voisinage ouvert $S' \subset S$ de $s$ tel que les images inverses de $z$ et $z'$ sur $S' \times_S X$ sont égales.*

On obtiendrait une condition équivalente en remplaçant $S'$ par un voisinage étale de $s$, c'est-à-dire en demandant l'existence d'un morphisme étale $S' \to S$ et l'existence d'un $s' \in S'(k)$ au-dessus de $s$ tel que $z$ et $z'$ coïncident sur $X' := S' \times_S X$.

En appliquant cette définition au morphisme $T^\star X^{\langle n+n' \rangle} \to X^{\langle n+n' \rangle}$ pour $X$ une courbe quasi-projective lisse sur $k$, nous allons démontrer le théorème 1.3.1 en montrant que les cycles $\tau^\star_{\mathscr{A}} \cdot \tau^\star_{\mathscr{A}'}$ et $\sum_{\mathscr{C} \in J(\mathscr{A}, \mathscr{A}')} \tau^\star_{\mathscr{C}}$ sur $T^\star X^{\langle n+n' \rangle}$ coïncident au voisinage de tout point $D \in X^{\langle n+n' \rangle}(k)$.

**1.3.9. Lemme.** *Soit $p : X' \to X$ un morphisme étale entre schémas séparés sur $k$. Le morphisme canonique $X' \simeq \Delta_{X'} \to X' \times_X X'$ est une immersion ouverte et fermée.*

Cela résulte par exemple de [ÉGA IV, 17.4.2].

**1.3.10. Lemme.** *Soit $p : X' \to X$ un morphisme étale entre courbes quasi-projectives lisses sur un corps algébriquement clos $k$. Notons $Z_p := X' \times_X X' - \Delta_{X'}$. Soit $I$ un ensemble fini de cardinal $n$. Pour tous $(i, j) \in I^2$ tels que $i \neq j$, notons $Z_{i,j} \subset X'^I$ l'image inverse de $Z_p$ par la projection $\pi_{i,j} : X'^I \to X'^2$. Notons $Z := \bigcup_{i \neq j} Z_{i,j}$, qui est fermé dans $X'^I$. Notons $\overline{Z} \subset X'^{\langle n \rangle}$ son image par le morphisme canonique $X'^I \to X'^{\langle n \rangle}$. Alors, les carrés du diagramme suivant sont cartésiens :*

$$\begin{array}{ccc} X'^{\mathscr{A}} - (X'^{\mathscr{A}} \cap Z) & \longrightarrow & X^{\mathscr{A}} \\ \downarrow & & \downarrow \\ X'^I - Z & \longrightarrow & X^I \\ \downarrow & & \downarrow \\ X'^{\langle n \rangle} - \overline{Z} & \longrightarrow & X^{\langle n \rangle} \end{array}.$$

*De plus, $X'^{\langle n\rangle} - \overline{Z}$ est précisément le lieu ouvert sur lequel le morphisme $X'^{\langle n\rangle} \to X^{\langle n\rangle}$ est étale.*

Dans le carré du haut, il est évident que les morphismes horizontaux sont étales. Il en résulte que le morphisme induit

$$(X'^{\mathscr{A}} - (X'^{\mathscr{A}} \cap Z) \to (X'^{I} - Z) \times_{X^I} X^{\mathscr{A}}$$

est étale ; par ailleurs, par construction de $Z$, ce morphisme induit une bijection sur les ensembles de $k$-points : c'est donc un isomorphisme.

Le fait que le lieu où $X'^{\langle n\rangle} \to X^{\langle n\rangle}$ soit étale soit exactement $X'^{\langle n\rangle} - \overline{Z}$ résulte de la description concrète des espaces tangents (cf. §A.2.3). Ceci étant acquis, il est évident que le morphisme canonique

$$X'^{I} - Z \to \left(X'^{\langle n\rangle} - \overline{Z}\right) \times_{X^{\langle n\rangle}} X^{I}$$

est étale, donc c'est un isomorphisme puisqu'il induit par ailleurs une bijection sur les ensembles de $k$-points.

**1.3.11. Lemme.** *Soit $p : X' \to X$ un morphisme étale entre courbes quasi-projectives lisses sur un corps algébriquement clos $k$. Supposons que $\mathscr{A}$ et $\mathscr{A}'$ sont des partitions de deux ensembles disjoints $I$ et $I'$ de cardinaux respectifs $n$ et $n'$. Soit $D' = \sum_{j\in J} \lambda_j x'_j$ un diviseur effectif de degré $n+n'$ sur $X'$ avec les $x'_j \in X'(k)$ distincts. On suppose que les points $x_j := p(x'_j) \in X(k)$ sont aussi distincts et on note $D := \sum_{j\in J} \lambda_j x_j$. Alors, les deux membres de l'égalité*

$$\tau^{\star}_{\mathscr{A}} \cdot \tau^{\star}_{\mathscr{A}'} = \sum_{\mathscr{C}\in J(\mathscr{A},\mathscr{A}')} \tau^{\star}_{\mathscr{C}}$$

*énoncée dans le théorème 1.3.1 pour la courbe $X'$ sont des cycles sur $T^{\star}X'^{\langle n+n'\rangle}$ qui coïncident au voisinage de $D' \in X^{\langle n+n'\rangle}(k)$ si et seulement si les deux cycles intervenant dans l'énoncé pour la courbe $X$ coïncident au voisinage de $D \in X^{\langle n+n'\rangle}(k)$.*

Avec les notations du lemme 1.3.10, on a $D' \in (X'^{\langle n+n'\rangle} - \overline{Z})(k)$, donc le morphisme $X'^{\langle n+n'\rangle} \to X^{\langle n+n'\rangle}$ est étale en $D'$. Au voisinage de $D'$, on peut donc identifier $T^{\star}X'^{\langle n+n'\rangle}$ et l'image inverse de $T^{\star}X^{\langle n+n'\rangle}$. Le lemme précédent implique aussi que pour toute partition $\mathscr{C}$ de $I \sqcup I'$, on a un isomorphisme :

$$X'^{\mathscr{C}} - (X'^{\mathscr{C}} \cap Z) \xrightarrow{\sim} (X'^{\langle n+n'\rangle} - \overline{Z}) \times_{X^{\langle n+n'\rangle}} X^{\mathscr{C}}.$$

On en déduit facilement que l'image inverse du cycle $\tau^{\star}_{\mathscr{C}} \in Z^{n+n'}(T^{\star}X^{\langle n+n'\rangle})$ sur $T^{\star}(X'^{\langle n+n'\rangle} - \overline{Z})$ coïncide avec la restriction du cycle $\tau^{\star}_{\mathscr{C}} \in Z^{n+n'}(T^{\star}X'^{\langle n+n'\rangle})$. Pour les mêmes raisons, cette assertion vaut aussi pour le produit $\tau^{\star}_{\mathscr{A}} \cdot \tau^{\star}_{\mathscr{A}'}$. On en déduit aussitôt la conclusion voulue.

**1.3.12.** Nous sommes maintenant en mesure de démontrer le théorème 1.3.1. Soit $X$ une courbe quasi-projective lisse sur un corps algébriquement clos $k$. Supposons que $\mathscr{A}$ et $\mathscr{A}'$ sont des partitions de deux ensembles disjoints $I$ et $I'$ de cardinaux respectifs $n$ et $n'$. Montrons que l'égalité de cycles sur $T^{\star}X^{\langle n+n'\rangle}$

$$\tau^{\star}_{\mathscr{A}} \cdot \tau^{\star}_{\mathscr{A}'} = \sum_{\mathscr{C}\in J(\mathscr{A},\mathscr{A}')} \tau^{\star}_{\mathscr{C}}$$

est vraie au voisinage de tout diviseur effectif $D$ de degré $n+n'$ sur $X$. On peut écrire $D = \sum_{j\in J} \lambda_j x_j$ avec des points $x_j \in X(k)$ distincts. Pour tout $j \in J$, choisissons arbitrairement des nombres distincts $z_j \in k$, un voisinage ouvert $U_j$ de $x_j$ dans $X$ et un morphisme étale $f_j : U_j \to \mathbf{A}^1_k$ envoyant $z_j$ sur $x_j$. Notons $X' := \coprod_{j\in J} U_j$ et $p\colon X' \to X$ le morphisme étale évident et $f\colon X' \to \mathbf{A}^1_k$ le morphisme correspondant aux morphismes $f_j$. On dispose d'un diagramme de morphismes étales :

$$\begin{array}{ccc} X' & \xrightarrow{p} & X \\ {\scriptstyle f}\downarrow & & \\ \mathbf{A}^1_k & & \end{array}$$

Notons $D'$ le diviseur $D' := \sum_{j\in J} \lambda_j \tilde{x}_j$ sur $X'$ où $\tilde{x}_j$ est le point de $U_j \subset X'$ correspondant à $x_j$. Notons aussi $D'' := \sum_{j\in J} \lambda_j z_j$ le diviseur sur $\mathbf{A}^1_k$. Le fait que les $(x_j)$, $(\tilde{x}_j)$, $(z_j)$ soient des familles de points distincts implique que l'on peut appliquer le lemme 1.3.11 pour les morphismes $f$ et $p$ avec le diviseur $D'$ sur $X'$. D'après la proposition 1.3.2, l'égalité voulue est satisfaite au voisinage de $D''$ pour la courbe $\mathbf{A}^1_k$. Par le lemme précédent appliqué à $f$, elle est satisfaite au voisinage de $D'$ pour $X'$. En l'appliquant à $p$, on obtient qu'elle est satisfaite au voisinage de $D$ pour la courbe $X$.

## 1.4. Formulaire.

**1.4.1. Proposition.** *Soit $X$ une courbe quasi-projective lisse sur un corps algébriquement clos $k$. Soit $(I_1,\dots,I_d)$ une famille d'ensembles finis disjoints. Soit $(\mathscr{A}_1,\dots,\mathscr{A}_d)$ une famille de partitions, où pour tout $i\in\{1,\dots,d\}$, $\mathscr{A}_i$ est une partition de $I_i$. On note $J(\mathscr{A}_1,\dots,\mathscr{A}_d)$ l'ensemble des partitions $\mathscr{C}$ de $J := I_1 \sqcup \dots \sqcup I_d$ qui pour tout $i\in\{1,\dots,d\}$ induisent la partition $\mathscr{A}_i$ et $I_i$, ce qui revient à demander que pour tout $C\in\mathscr{C}$ et tout $i\in\{1,\dots,d\}$, $C\cap I_i \in \mathscr{A}_i \cup \{\varnothing\}$. Alors, on a l'égalité*

$$\prod_{i=1}^{d} \tau^\star_{\mathscr{A}_i} = \sum_{\mathscr{C}\in J(\mathscr{A}_1,\dots,\mathscr{A}_d)} \tau^\star_{\mathscr{C}} \in Z^{\#J}(T^\star X^{\langle \#J\rangle}).$$

Cette proposition se démontre par récurrence évidente sur $d$, le cas essentiel $d=2$ étant le théorème 1.3.1.

**1.4.2. Définition.** *Si $\mathscr{A}$ est une partition d'un ensemble fini de cardinal $n$, que l'on note $\lambda_1 \geq \lambda_2 \geq \dots \geq \lambda_l \geq 1$ avec $n = \lambda_1 + \dots + \lambda_l$ les cardinaux des parties de la partition $\mathscr{A}$ et que l'on note $\underline{e} = (e_1,\dots,e_n)$ où $e_i := \#\{j\in\{1,\dots,l\}, \lambda_j = i\} = \#\{A\in\mathscr{A}, \#A = i\}$ alors on a défini deux cycles $\tau^\star_{\mathscr{A}}$ et $\tau^\star_{\underline{e}}$ qui vérifient la relation $\tau^\star_{\mathscr{A}} = \underline{e}!\tau^\star_{\underline{e}}$ où $\underline{e}! = \prod_{i=1}^{n} e_i!$. Le cycle algébrique $\tau^\star_{\underline{e}}$ sera aussi noté $\tau^\star_{n=\lambda_1+\dots+\lambda_l}$.*

Par exemple, $\tau^\star_{4=2+1+1} = \tau^\star_{(2,1)} = \frac{1}{2}\tau^\star_{\mathscr{A}}$ où $\mathscr{A} = \{\{1,2\},\{3\},\{4\}\}$.

On prendra garde à ne pas confondre $\tau^\star_{(n)} = \tau^\star_{n=1+1+\dots+1}$ et $\tau^\star_{n=n} = \tau^\star_{\underline{e}}$ où $\underline{e} = (0,\dots,0,1)$ vérifie $e_n = 1$.

**1.4.3. Définition.** *Supposons que $\underline{\lambda} = (\lambda_1,\dots,\lambda_l)$ soit une partition de l'entier $n$ et que $\underline{\mu} = (\mu_1,\dots,\mu_m)$ soit une composition de $n$ (c'est-à-dire que les entiers $\mu_i$ sont non nuls et tels que $|\underline{\mu}| := \mu_1 + \dots + \mu_m = n$, mais ne sont pas nécessairement en ordre décroissant). On note $m_{\underline{\lambda},\underline{\mu}} \in \mathbf{N}$ le cardinal de l'ensemble $\mathscr{X}_{\underline{\lambda},\underline{\mu}}$ des parties $Z \subset \{1,\dots,l\}\times\{1,\dots,m\}$ telles que*

- *pour tout $i\in\{1,\dots,l\}$, $\#\{j\in\{1,\dots,m\}, (i,j)\in Z\} = \lambda_i$ ;*
- *pour tout $j\in\{1,\dots,m\}$, $\#\{i\in\{1,\dots,l\}, (i,j)\in Z\} = \mu_j$.*

La proposition suivante généralise une formule obtenue par Gérard Laumon en caractéristique nulle [Laumon 1987a, p. 330-331] :

**1.4.4. Proposition.** *Soit $X$ une courbe quasi-projective lisse sur un corps algébriquement clos $k$. Soit $\Lambda$ un corps fini de caractéristique inversible dans $k$. Soit $n\in\mathbf{N}$. Soit $\underline{\mu} = (\mu_1,\dots,\mu_m)$ une composition* *de l'entier n. Notons* $\pi_{\underline{\mu}} \colon X^{\langle \underline{\mu} \rangle} = \prod_{j=1}^{m} X^{\langle \mu_i \rangle} \to X^{\langle n \rangle}$ *le morphisme canonique. Alors,*

$$CC(X^{\langle n\rangle}, \pi_{\underline{\mu},\star}\Lambda) = (-1)^n \tau^{\star}_{(\mu_1)} \cdot \cdots \cdot \tau^{\star}_{(\mu_m)} = (-1)^n \sum_{\underline{\lambda} \text{ partition de } n} m_{\underline{\lambda},\underline{\mu}} \tau^{\star}_{n=\lambda_1+\cdots+\lambda_l}.$$

*En outre, le support singulier* $SS(X^{\langle n\rangle}, \pi_{\underline{\mu},\star}\Lambda)$ *est exactement le support du cycle caractéristique.*

On peut plonger $X$ comme ouvert dans $\overline{X}$ projective lisse. Quitte à remplacer $X$ par $\overline{X}$, on peut supposer que $X$ est projective lisse. On remarque que pour tout $j \in \{1, \dots, m\}$, on a :

$$CC(X^{\langle \mu_j\rangle}, \Lambda) = (-1)^{\mu_j} \tau^{\star}_{(\mu_j)}$$

et le support singulier est exactement cette section nulle. Comme $\pi_{\underline{\mu},\star}\Lambda \simeq \vee_{j=1}^{m} \Lambda_{|X^{\langle \mu_j\rangle}}$, la proposition **1.2.4** montre que

$$CC(X^{\langle n\rangle}, \pi_{\underline{\mu},\star}\Lambda) = \prod_{j=1}^{m} \left( (-1)^{\mu_j} \tau^{\star}_{(\mu_j)} \right) = (-1)^n \prod_{j=1}^{m} \tau^{\star}_{(\mu_j)}.$$

Pour montrer la dernière égalité voulue, choisissons $\{I_1, \dots, I_m\}$ une partition de $\{1, \dots, n\}$ telle que $\#I_j = \mu_j$ pour tout $j \in \{1, \dots, m\}$. Notons $\mathscr{A}_j$ la partition de $I_j$ en $\mu_j$ singletons. On a alors : $\tau^{\star}_{\mathscr{A}_j} = \mu_j! \tau^{\star}_{(\mu_j)}$ pour tout $j \in \{1, \dots, m\}$ et la proposition **1.4.1** implique que

$$\left( \prod_{j=1}^{m} \mu_j! \right) \cdot \prod_{j=1}^{m} \tau^{\star}_{(\mu_j)} = \prod_{j=1}^{m} \tau^{\star}_{\mathscr{A}_j} = \sum_{\mathscr{C} \in J(\mathscr{A}_1, \dots, A_m)} \tau^{\star}_{\mathscr{C}}.$$

Soit $\underline{\lambda} = (\lambda_1, \dots, \lambda_l)$ une partition de $n$. Notons $P_{\underline{\lambda}}$ l'ensemble des partitions $\mathscr{C} \in J(\mathscr{A}_1, \dots, A_m)$ dont la partition de l'entier $n$ associée soit $\underline{\lambda}$. Notons $\underline{e} = (e_1, \dots, e_n)$ où $e_s \coloneqq \#\{i \in \{1, \dots, l\}, \lambda_i = s\}$ pour tout $s \in \{1, \dots, n\}$ de sorte que pour tout $\mathscr{C} \in P_{\lambda}$, on ait :

$$\tau^{\star}_{\mathscr{C}} = \underline{e}! \cdot \tau^{\star}_{\underline{e}} = \underline{e}! \cdot \tau^{\star}_{n=\lambda_1+\cdots+\lambda_l}.$$

Il s'agit donc de montrer l'égalité suivante :

$$\#P_{\underline{\lambda}} \cdot \underline{e}! = \left( \prod_{j=1}^{m} \mu_j! \right) \cdot m_{\underline{\lambda},\underline{\mu}}.$$

On introduit l'ensemble $\widetilde{P}_{\underline{\lambda}}$ formé des $m$-uplets $(C_1, \dots, C_l)$ de parties disjointes de $\{1, \dots, n\}$ telles que $\mathscr{C} \coloneqq \{C_1, \dots, C_l\} \in J(\mathscr{A}_1, \dots, A_m)$ et que pour tout $j \in \{1, \dots, m\}$, $\#C_j = \mu_j$.

L'application évidente $\widetilde{P}_{\underline{\lambda}} \to P_{\underline{\lambda}}$ permet de montrer que $\#\widetilde{P}_{\underline{\lambda}} = \underline{e}! \cdot \#P_{\underline{\lambda}}$.

Considérons l'application

$$\widetilde{P}_{\underline{\lambda}} \to \mathscr{X}_{\underline{\lambda},\underline{\mu}}$$

qui à $(C_1, \dots, C_l)$ associe $Z \coloneqq \{(i,j) \in \{1, \dots, l\} \times \{1, \dots, m\}, I_j \cap C_i \neq \varnothing\}$.

Avec ces notations, on remarque que la condition $\mathscr{C} = \{C_1, \dots, C_l\} \in J(\mathscr{A}_1, \dots, \mathscr{A}_m)$ implique que $I_j \cap C_i$ est soit l'ensemble vide soit un singleton que l'on notera $\{x_{i,j}\}$ pour tout $(i,j) \in Z$. On obtient ainsi une bijection $Z \xrightarrow{\sim} \{1, \dots, n\}$ qui à $(i,j)$ associe $x_{i,j}$. Pour tout $j \in \{1, \dots, m\}$, cette bijection a la propriété d'induire une bijection $Z_j \xrightarrow{\sim} I_j$ où $Z_j = \{(i,k) \in Z, k = j\}$.

Réciproquement, si on fixe $Z \in \mathscr{X}_{\underline{\lambda},\underline{\mu}}$, il existe $\prod_{j=1}^{m} \mu_j!$ manières de choisir des bijections $Z_j \xrightarrow{\sim} I_j$ pour tout $j \in \{1, \dots, m\}$ et pour chacun de ces choix, on a une bijection $x \colon Z = \coprod_j Z_j \xrightarrow{\sim} \{1, \dots, n\} = \coprod_j I_j$ qui permet d'obtenir un antécédent $(C_1, \dots, C_l) \in \widetilde{P}_{\underline{\lambda}}$ de $Z$ en notant pour tout $i \in \{1, \dots, l\}$, $C_i$ l'ensemble des $x_{i,j}$ pour $j \in \{1, \dots, m\}$ tel que $(i,j) \in Z$.

On a montré l'identité suivante :

$$\#\widetilde{P}_{\underline{\lambda}} = \left(\prod_{j=1}^{m} \mu_j!\right) \cdot m_{\underline{\lambda},\underline{\mu}}.$$

**1.4.5. Remarque.** Si on note $\underline{\lambda}^{\mathrm{T}}$ la partition conjuguée d'une partition $\underline{\lambda}$, c'est-à-dire celle dont le diagramme de Young est le transposé de celui de $\underline{\lambda}$, et que l'on ordonne les partitions d'un entier $n$ par l'ordre lexicographique inverse, alors on montre facilement que la matrice $(m_{\underline{\lambda}^{\mathrm{T}},\underline{\mu}})_{\underline{\lambda},\underline{\mu}}$, où $\underline{\lambda}$ et $\underline{\mu}$ parcourent les partitions de $n$, est triangulaire supérieure avec des 1 sur la diagonale. En utilisant la proposition précédente, on obtient que tous les cycles $\tau^\star_{\underline{\lambda}}$ peuvent être obtenus comme cycles caractéristiques d'objets de $\mathsf{D}^{\mathsf{b}}_{\mathsf{c}}(X^{\langle n\rangle}, \Lambda)$.

**1.4.6. Proposition.** *Soit $X$ une courbe quasi-projective lisse sur un corps algébriquement clos $k$. Soit $\Lambda$ un corps fini de caractéristique inversible dans $k$. Soit $\mathscr{A}$ une partition d'un ensemble fini $\boldsymbol{I}$ de cardinal $n$. Notons $\pi^{\mathscr{A}} : X^{\mathscr{A}} \to X^{\boldsymbol{I}} \simeq X^n \to X^{\langle n\rangle}$ le morphisme composé. Notons $\underline{\lambda} = (\lambda_1, \dots, \lambda_l)$ la partition de $n$ associée à $\mathscr{A}$ que l'on écrit $\mathscr{A} = \{A_1, \dots, A_l\}$ avec $\#A_i = \lambda_i$ pour tout $i \in \{1, \dots, l\}$.*

*Si pour tout $i \in \{1, \dots, l\}$, l'entier $\lambda_i$ est inversible dans $k$, alors on a l'égalité suivante entre cycles sur $T^\star X^{\langle n\rangle}$ :*

$$CC(X^{\langle n\rangle}, \pi^{\mathscr{A}}_\star \Lambda) = (-1)^l \sum_{\mathscr{B} \geq \mathscr{A}} \tau^\star_{\mathscr{B}}$$

*où l'ensemble d'indices est l'ensemble des partitions $\mathscr{B}$ de $\boldsymbol{I}$ telles que pour tout $\boldsymbol{B} \in \mathscr{B}$, $\boldsymbol{B}$ soit une réunion de parties appartenant à $\mathscr{A}$. De plus, le support singulier de $\pi^{\mathscr{A}}_\star \Lambda$ est exactement le support du cycle caractéristique.*

*Démonstration.* Pour tout $i \in \{1, \dots, l\}$, le fait que $\lambda_i$ soit inversible dans $k$ a pour conséquence que le morphisme évident $\iota_i : X \xrightarrow{\Delta} X^{\lambda_i} \to X^{\langle \lambda_i \rangle}$ est une immersion fermée (cf. corollaire 1.1.6) et

$$CC(X^{\langle \lambda_i\rangle}, \iota_{i\star}\Lambda) = -\tau^\star_{\lambda_i = \lambda_i}.$$

En utilisant la proposition 1.2.4 (dans le cas $X$ projective lisse, cas auquel on se ramène), il vient alors :

$$CC(X^{\langle n\rangle}, \pi^{\mathscr{A}}_\star \Lambda) = CC(X^{\langle n\rangle}, \vee_{i=1}^{l} \iota_{i\star}\Lambda) = \prod_{i=1}^{l} \left(-\tau^\star_{\lambda_i=\lambda_i}\right) = (-1)^l \cdot \prod_{i=1}^{l} \tau^\star_{\lambda_i=\lambda_i}$$

Pour tout $i \in \{1, \dots, l\}$, notons $\mathscr{A}_i := \{A_i\}$ la partition évidente de l'ensemble $A_i$ en une seule partie à $\lambda_i$ éléments. D'après la proposition 1.4.1, on a :

$$\prod_{i=1}^{l} \tau^\star_{\lambda_i=\lambda_i} = \prod_{i=1}^{l} \tau_{\mathscr{A}_i} = \sum_{\mathscr{B} \in J(\mathscr{A}_1, \dots, \mathscr{A}_l)} \tau^\star_{\mathscr{B}}$$

On peut conclure puisque l'ensemble d'indices $J(\mathscr{A}_1, \dots, \mathscr{A}_l)$ est exactement celui qui a été décrit dans l'énoncé. □

**1.4.7. Remarque.** En caractéristique zéro, grâce à cette proposition, on construit facilement par récurrence sur $\#\mathscr{A}$ des objets de $\mathsf{D}^{\mathsf{b}}_{\mathsf{c}}(X^{\langle n\rangle}, \Lambda)$ dont l'image directe par $X^{\langle n\rangle} \to X^n$ ait pour cycle caractéristique le cycle $\tau^\star_{\mathscr{A}}$.

**1.4.8. Proposition.** *Soit $X$ une courbe quasi-projective lisse sur un corps algébriquement clos $k$. Soit $\Lambda$ un corps fini de caractéristique inversible dans $k$. Soit $\mathscr{F}$ un faisceau de $\Lambda$-modules localement constant de rang $r$ sur $X$. Alors, si on note $S_{\mathscr{F}} = \sum_{n=0}^{\infty} CC(X^{\langle n\rangle}, \mathscr{F}^{\langle n\rangle})$, on a :*

$$S_{\mathscr{F}} = S_{\Lambda^r} = (S_\Lambda)^r = \sum_{n=0}^{\infty} (-1)^n \sum_{\underline{e}, \|\underline{e}\|=n} \left(\prod_{i=1}^{n} \binom{r}{i}^{e_i}\right) \cdot \tau^\star_{\underline{e}}$$

*et pour tout $n \in \mathbf{N}$, le support singulier de $\mathscr{F}^{\langle n \rangle}$ est exactement le support son cycle caractéristique.*

Commençons par le cas $\mathscr{F} = \Lambda$. On a alors pour tout $n \in \mathbf{N}$ un isomorphisme canonique $\Lambda^{\langle n \rangle} \simeq \Lambda$ de faisceaux sur $X^{\langle n \rangle}$, donc :

$$S_\Lambda = \sum_{n=0}^{\infty} CC(X^{\langle n \rangle}, \Lambda) = \sum_{n=0}^{\infty} (-1)^n \tau^\star_{(n)}.$$

Pour montrer l'égalité $S_{\Lambda^r} = (S_\Lambda)^r$, on peut supposer que $X$ est projective lisse et appliquer la proposition **1.2.4**. On a alors :

$$\begin{aligned} S_{\Lambda^r} &= (S_\Lambda)^r = \prod_{j=1}^{r} \left( \sum_{\mu_j=0}^{\infty} (-1)^{\mu_j} \tau^\star_{(\mu_j)} \right) \\ &= \sum_{(\mu_1,\dots,\mu_r)\in\mathbf{N}^r} (-1)^{\sum_{j=1}^r \mu_j} \prod_{j=1}^{r} \tau^\star_{(\mu_j)} \\ &= \sum_{n=0}^{\infty} (-1)^n \sum_{\substack{(\mu_1,\dots,\mu_r)\in\mathbf{N}^r \\ |\underline{\mu}|=n}} \prod_{j=1}^{r} \tau^\star_{(\mu_j)} \end{aligned}$$

D'après la proposition **1.4.4**, on obtient :

$$\begin{aligned} S_{\Lambda^r} &= \sum_{n=0}^{\infty} (-1)^n \sum_{\substack{(\mu_1,\dots,\mu_r)\in\mathbf{N}^r \\ |\underline{\mu}|=n}} \sum_{\underline{\lambda} \text{ partition de } n} m_{\underline{\lambda},\underline{\mu}} \tau^\star_{n=\lambda_1+\dots+\lambda_l} \\ &= \sum_{n=0}^{\infty} (-1)^n \sum_{\underline{\lambda} \text{ partition de } n} \left( \sum_{\substack{(\mu_1,\dots,\mu_r)\in\mathbf{N}^r \\ |\underline{\mu}|=n}} m_{\underline{\lambda},\underline{\mu}} \right) \tau^\star_{n=\lambda_1+\dots+\lambda_l} \end{aligned}$$

Pour $\underline{\lambda} = (\lambda_1, \dots, \lambda_l)$ fixé, la somme des $m_{\underline{\lambda},\underline{\mu}}$ pour $(\mu_1, \dots, \mu_r) \in \mathbf{N}^r$ tels que $|\underline{\mu}| = n = |\underline{\lambda}|$ est la somme des cardinaux des ensembles disjoints $\mathscr{X}_{\underline{\lambda},\underline{\mu}}$ (cf. définition **1.4.3**) dont la réunion est l'ensemble $\mathscr{X}_{\underline{\lambda},d}$ des parties $Z \subset \{1, \dots, l\} \times \{1, \dots, r\}$ telles que pour tout $i \in \{1, \dots, l\}$, $\#\{j \in \{1, \dots, r\}, (i, j) \in Z\} = \lambda_i$. Ainsi, il est évident que :

$$\sum_{\substack{(\mu_1,\dots,\mu_r)\in\mathbf{N}^r \\ |\underline{\mu}|=n}} m_{\underline{\lambda},\underline{\mu}} = \#\mathscr{X}_{\underline{\lambda},r} = \prod_{i=1}^{l} \binom{r}{\lambda_i}.$$

Par conséquent,

$$\begin{aligned} S_{\Lambda^r} &= \sum_{n=0}^{\infty} (-1)^n \sum_{\underline{\lambda} \text{ partition de } n} \left( \prod_{i=1}^{l} \binom{r}{\lambda_i} \right) \tau^\star_{n=\lambda_1+\dots+\lambda_l} \\ &= \sum_{n=0}^{\infty} (-1)^n \sum_{\underline{e}, \|\underline{e}\|=n} \left( \prod_{i=1}^{n} \binom{r}{i}^{e_i} \right) \tau^\star_{\underline{e}} \end{aligned}$$

Ceci démontre la formule voulue dans le cas des faisceaux constants. Dans le cas général où $\mathscr{F}$ est un faisceau localement constant de rang $r$, on raisonne comme en **1.3.12** : il suffit de montrer que la formule pour $CC(X^{\langle n \rangle}, \mathscr{F}^{\langle n \rangle})$ (et son support singulier) est vraie au voisinage de tout point $D \in X^{\langle n \rangle}(k)$. On peut écrire $D = \sum_{j\in J} \nu_j x_j$ avec des points $x_j \in X(k)$ distincts. Pour chaque $j \in J$, on choisit $f_j : V_j \to X$ étale avec $v_j \in V_j(k)$ au-dessus de $x_j$ tel que $f_j^\star \mathscr{F}$ soit isomorphe à $\Lambda^r$. On peut noter $\widetilde{X} := \sqcup V_j$ et $\widetilde{D} := \sum_{j\in J} \lambda_j v_j$. De même que dans la situation du lemme 1.3.11, en utilisant le morphisme étale évident $f\colon \widetilde{X} \to X$, la formule affirmée dans la proposition pour $CC(X^{\langle n\rangle}, \mathscr{F}^{\langle n\rangle})$ est vraie au voisinage de $D$ si et seulement si la formule pour $CC(\widetilde{X}^{\langle n\rangle}, (f^\star\mathscr{F})^{\langle n\rangle})$ est vraie au voisinage de $\widetilde{D}$, ce qui est bien le cas puisque $f^\star\mathscr{F}$ est isomorphe à $\Lambda^r$.

**1.4.9. Définition.** *Si $\Delta = \sum_{s\in S} m_s \cdot s$ est un diviseur effectif sur une courbe quasi-projective lisse $X$ sur un corps algébriquement clos $k$ (avec $S$ une partie finie de $X(k)$) et que $\Delta' = \sum_{s\in S} a_s \cdot s$ est un diviseur (éventuellement non effectif), on notera :*

$$\binom{\Delta'}{\Delta} := \prod_{s\in S} \binom{a_s}{m_s}.$$

(Conformément à l'usage, on étend la définition des coefficients binomiaux $\binom{a}{m}$ au cas où $a \in \mathbf{Z}$ (et $m \in \mathbf{N}$) par la formule $\frac{\prod_{i=0}^{m-1}(a-i)}{m!} \in \mathbf{Z}$. On a alors $\binom{-a}{m} = \binom{a+m-1}{m}$.)

**1.4.10. Proposition.** *Soit $X$ une courbe quasi-projective lisse sur un corps algébriquement clos $k$. Soit $\Lambda$ un corps fini de caractéristique inversible dans $k$. Soit $U$ un ouvert dense de $X$. Notons $S := X(k) - U(k)$. Soit $\mathscr{P}$ un faisceau constructible de $\Lambda$-modules tel que $\mathscr{P}_{|U} \simeq 0$. On note $\Delta_{\mathscr{P}} := \sum_{s\in S} r_s \cdot s$ où $r_s := \dim_\Lambda \mathscr{P}_s$. On a alors :*

$$S_{\mathscr{P}} = \prod_{s\in S}(1-\tau_s^\star)^{-r_s} = \sum_{\substack{\Delta=\sum_{s\in S} m_s\cdot s \\ \text{diviseur effectif}}} (-1)^{|\Delta|}\binom{-\Delta_{\mathscr{P}}}{\Delta}\tau_\Delta^\star$$

$$S_{\mathscr{P}[-1]} = \prod_{s\in S}(1-\tau_s^\star)^{r_s} = \sum_{0\le\Delta\le\Delta_{\mathscr{P}}} (-1)^{|\Delta|}\binom{\Delta_{\mathscr{P}}}{\Delta}\tau_\Delta^\star$$

*En outre, pour tout $n \in \mathbf{N}$, le support singulier des faisceaux $\mathscr{P}^{\langle n\rangle}$ et $\mathscr{P}[-1]^{\langle n\rangle}$ est exactement le support de leurs cycles caractéristiques respectifs.*

Commençons par le cas où $S = \{s\}$ est un singleton. Pour tout diviseur effectif $D$ de degré $n$, on notera $i_D\colon \operatorname{Spec}(k) \to X^{\langle n\rangle}$ l'immersion fermée associée. Si $F$ est un complexe borné de $\Lambda$-espaces vectoriels, il est évident que l'on a un isomorphisme canonique dans $\mathrm{D^b}(X^{\langle n\rangle}, \Lambda)$ :

$$(i_{s,\star}F)^{\langle n\rangle} \simeq i_{ns,\star}(F^{\langle n\rangle})$$

où les puissances symétriques $F^{\langle n\rangle} = \mathsf{Y}^n(F)$ (cf. A.1.1, A.1.4) sont calculées dans la catégorie dérivée des $F$-espaces vectoriels. En outre, on dispose de la formule de décalage de Quillen-Illusie [Illusie 1971-1972, I.4.3.2.1] :

$$F[-1]^{\langle n\rangle} \simeq (L\wedge^n F)[-n].$$

On en déduit aussitôt que si on note $\mathscr{P} := i_{s,\star}F$ avec $F$ un $\Lambda$-espace vectoriel de dimension finie $r_s$, on a un isomorphisme $\mathscr{P}[-1]^{\langle n\rangle} \simeq (i_{ns,\star}\wedge^n F)[-n]$ pour tout $n \in \mathbf{N}$. La dimension de $\wedge^n F$ étant $\binom{r_s}{n}$, on obtient :

$$S_{\mathscr{P}[-1]} = \sum_{n=0}^{r_s}(-1)^n\binom{r_s}{n}\tau_{ns}^\star = (1-\tau_s^\star)^{r_s}$$

et le support singulier des faisceaux $\mathscr{P}[-1]^{\langle n\rangle}$ est exactement le support du cycle caractéristique.

Pour tout $n \in \mathbf{N}$, on a aussi $\mathscr{P}^{\langle n\rangle} \simeq i_{ns,\star}(\Gamma^n F)$ dont le support singulier est $|\tau_{ns}^\star|$ sauf si $F = 0$ et $n \ge 1$, auquel cas il est vide. Quitte à se ramener cas où $X$ est projective lisse, on peut utiliser la proposition 1.2.13 et obtenir le résultat voulu :

$$S_{\mathscr{P}} = \frac{1}{S_{\mathscr{P}[-1]}} = (1-\tau_s^\star)^{-r_s} = \sum_{n=0}^{\infty}(-1)^n\binom{-r_s}{n}\tau_{ns}^\star.$$

Dans le cas général, le faisceau $\mathscr{P}$ s'écrit $\mathscr{P} = \bigoplus_{s\in S}\mathscr{P}_s$ où $\mathscr{P}_s$ est un faisceau gratte-ciel supporté par $s$. Comme on peut facilement supposer que $X$ est projective lisse, on peut conclure en utilisant la proposition 1.2.10.

**1.4.11. Proposition.** *Soit $X$ une courbe quasi-projective lisse sur un corps algébriquement clos $k$. Soit $\Lambda$ un corps fini de caractéristique inversible dans $k$. Soit $U$ un ouvert dense de $X$. Notons $S := X(k)-U(k)$ et $j : U \to X$ l'immersion ouverte. Soit $\mathscr{F}$ un faisceau constructible de $\Lambda$-modules sur $X$. On suppose que $j_\star j^\star F$ est un faisceau localement constant de rang $r$ et que $\mathscr{F} \to j_\star j^\star\mathscr{F}$ est injectif. Notons $\mathscr{P}$ le conoyau de cette injection, de sorte que l'on dispose d'une suite exacte courte $0 \to \mathscr{F} \to j_\star j^\star\mathscr{F} \to \mathscr{P} \to 0$.*

*Notons $\Delta_{\mathscr{F}} := \Delta_{\mathscr{P}} = \sum_{s\in S} a_s \cdot s$ où $a_s := \dim_\Lambda \mathscr{P}_s = r - \dim_\Lambda \mathscr{F}_s$. Alors,*

$$S_{\mathscr{F}} = S_{\Lambda^r} \cdot S_{\mathscr{P}[-1]}.$$

*Autrement dit, pour tout $n \in \mathbf{N}$, on a :*

$$CC(X^{\langle n\rangle}, \mathscr{F}^{\langle n\rangle}) = (-1)^n \sum_{\substack{\Delta,\underline{e},|\Delta|+\|\underline{e}\|=n\\ 0\le\Delta\le\Delta_{\mathscr{F}}}} \binom{\Delta_{\mathscr{F}}}{\Delta}\left(\prod_{i=1}^{n}\binom{r}{i}^{e_i}\right)\tau_{\Delta,\underline{e}}^\star$$

*et le support singulier $SS(X^{\langle n\rangle}, \mathscr{F}^{\langle n\rangle})$ est exactement le support du cycle caractéristique.*

Notons $\mathscr{G} := j_\star j^\star\mathscr{F}$, qui est un faisceau localement constant auquel on peut appliquer la proposition 1.4.8, ce qui permet de calculer la série $S_{\mathscr{G}} = S_{\Lambda^r}$. La proposition 1.4.10 donne une formule pour $S_{\mathscr{P}[-1]}$. On dispose aussi du triangle distingué :

$$\mathscr{P}[-1] \to \mathscr{F} \to \mathscr{G} \to \mathscr{P}.$$

Pour conclure, il suffit de montrer $S_{\mathscr{F}} = S_{\mathscr{G}} \cdot S_{\mathscr{P}[-1]}$. Dans le cas où $X$ est projective, cela résulte directement de la proposition 1.2.11. Cependant, la propriété de multiplicativité énoncée dans cette proposition s'applique ici même si $X$ est seulement quasi-projective puisque le faisceau $\mathscr{P}$ est à support ponctuel : en effet, on est ramené à montrer la formule de la proposition 1.2.10 dans le cas où un des deux objets considérés est un faisceau à support ponctuel. On peut alors conclure en utilisant la remarque 1.2.5.

## 1.5. Faisceaux modérés.

**1.5.1. Théorème.** *Soit $X$ une courbe quasi-projective lisse sur un corps algébriquement clos $k$. Soit $\Lambda$ un corps fini de caractéristique inversible dans $k$. Soit $U$ un ouvert dense de $X$. Notons $S := X(k) - U(k)$ et $j : U \to X$ l'immersion ouverte. Soit $\mathscr{F}$ un faisceau constructible sur $X$. On suppose que $j^\star\mathscr{F}$ est localement constant de rang $r$, que $\mathscr{F}$ est modéré le long de $S$ et que $\mathscr{F} \to j_\star j^\star\mathscr{F}$ est injectif. On note $\Delta_{\mathscr{F}} := \sum_{s\in S} a_s \cdot s$ où $a_s := r - \dim_\Lambda \mathscr{F}_s$. Alors,*

$$S_{\mathscr{F}} = S_{\Lambda^r} \cdot \prod_{s\in S}(1-\tau_s^\star)^{a_s}.$$

*Autrement dit, pour tout $n \in \mathbf{N}$, on a la formule :*

$$CC(X^{\langle n\rangle}, \mathscr{F}^{\langle n\rangle}) = (-1)^n \sum_{\substack{\Delta,\underline{e},|\Delta|+\|\underline{e}\|=n\\ 0\le\Delta\le\Delta_{\mathscr{F}}}} \binom{\Delta_{\mathscr{F}}}{\Delta}\left(\prod_{i=1}^{n}\binom{r}{i}^{e_i}\right)\tau_{\Delta,\underline{e}}^\star$$

*et le support singulier $SS(X^{\langle n\rangle}, \mathscr{F}^{\langle n\rangle})$ est exactement le support du cycle caractéristique.*

La démonstration du théorème occupe le reste de cette section. Nous allons commencer par établir le résultat dans certaines cas particuliers de faisceaux de rang 1 sur la droite projective :

**1.5.2. Proposition.** *Soit $k$ un corps algébriquement clos $k$. Soit $\Lambda$ un corps fini de caractéristique inversible dans $k$. Soit $\mathscr{F}$ un faisceau localement constant de $\Lambda$-modules de rang 1 sur $\mathbf{G}_{m_k}$. On suppose qu'il existe un entier $N$ inversible dans $k$ tel que $[N]^\star\mathscr{F} \simeq \Lambda$ où $[N]\colon \mathbf{G}_{m_k} \to \mathbf{G}_{m_k}$ est l'élévation à la puissance $N$. Notons $j\colon \mathbf{G}_{m_k} \to \mathbf{P}^1_k$ l'immersion ouverte canonique. Alors,*

$$S_{j_!\mathscr{F}} = S_{j_!\Lambda} = S_\Lambda \cdot (1 - \tau^\star_{[0]}) \cdot (1 - \tau^\star_{[\infty]}).$$

*En outre, pour tout $n \in \mathbf{N}$, le support singulier de $(j_!\mathscr{F})^{\langle n\rangle}$ est exactement le support de son cycle caractéristique.*

**1.5.3. Lemme.** *Sous les hypothèses de la proposition 1.5.2, si on note $\varphi\colon \mathbf{G}^2_{m_k} \to \mathbf{G}^2_{m_k}$ le morphisme $(x, y) \longmapsto (1, xy)$, alors on dispose d'un isomorphisme canonique de faisceaux sur $\mathbf{G}^2_{m_k}$ :*

$$\varphi^\star(\mathscr{F} \boxtimes \mathscr{F}) \simeq \mathscr{F} \boxtimes \mathscr{F}.$$

*(Cet isomorphisme est unique si on impose qu'il coïncide avec l'identité au point $(1, 1) \in (k^\times)^2$.)*

*Autrement dit, si on choisit un isomorphisme $\mathscr{F}_1 \simeq \Lambda$ et que l'on note $\mu\colon \mathbf{G}_{m_k}^{\,2} \to \mathbf{G}_{m_k}$ la multiplication, on a un isomorphisme canonique de faisceaux sur $\mathbf{G}^2_{m_k}$ :*

$$\mathscr{F} \boxtimes \mathscr{F} \simeq \mu^\star\mathscr{F}.$$

*Plus généralement, si $n \geq 0$, le faisceau $\mathscr{F}^{\boxtimes n}$ sur $\mathbf{G}^n_{m_k}$ s'identifie à l'image inverse de $\mathscr{F}$ par le morphisme de multiplication $\mathbf{G}^n_{m_k} \to \mathbf{G}_{m_k}$.*

*Démonstration.* Notons $\Gamma \simeq \mu_N(k)$ le groupe de Galois du revêtement étale $\mathbf{G}_{m_k} \xrightarrow{[N]} \mathbf{G}_{m_k}$. Le revêtement produit de $\mathbf{G}^2_{m_k}$ de groupe de Galois $\Gamma \times \Gamma$ est le schéma obtenu en extrayant des racines $n$-ièmes $X^{\frac{1}{n}}$ et $Y^{\frac{1}{n}}$ à $X$ et à $Y$ dans $\mathbf{G}^2_{m_k} = \mathrm{Spec}(k[X, Y, X^{-1}, Y^{-1}])$. La fonction $X^{\frac{1}{n}}Y^{\frac{1}{n}}$ est fixée par le sous-groupe $H$ de $\Gamma \times \Gamma$ formé des éléments de la forme $(\gamma, \gamma^{-1})$ pour $\gamma \in \Gamma$.

Le faisceau $\mathscr{F}$ est déterminé par l'action $\rho\colon \Gamma \to \mathrm{GL}_1(\Lambda) = \Lambda^\times$ sur le $\Lambda$-espace vectoriel de dimension un $\mathscr{F}_1$. Le faisceau $\mathscr{F} \boxtimes \mathscr{F}$ est déterminé par l'action produit $\Gamma \times \Gamma \to \Lambda^\times$ qui se factorise évidemment par le sous-groupe $H$. Ceci montre que $\mathscr{F} \boxtimes \mathscr{F}$ est trivialisé par le revêtement de $\mathbf{G}^2_{m_k} = \mathrm{Spec}(k[X, Y, X^{-1}, Y^{-1}])$ correspondant à ce sous-groupe, c'est-à-dire $W \coloneqq \mathrm{Spec}(k[X, Y, X^{-1}, Y^{-1}, U]/(U^n - XY)$.

Le faisceau $\mathscr{F} \boxtimes \mathscr{F}$ est obtenu par descente du faisceau constant de valeur $\mathscr{F}_1 \otimes \mathscr{F}_1$ sur $W$ relativement à l'action évidente de $(\Gamma \times \Gamma)/H \xrightarrow{\sim} \Gamma$. Comme le morphisme d'anneaux $k[X, Y, X^{-1}, Y^{-1}] \to k[X, Y, X^{-1}, Y^{-1}]$ correspondant à $\varphi$ envoie $XY$ sur $1 \cdot (XY) = XY$, et que $W$ est obtenu en extrayant une racine $n$-ième de $XY$ dans $\mathbf{G}^2_{m_k}$, on en déduit un diagramme cartésien évident compatible aux actions de $\Gamma$ :

$$\begin{array}{ccc} W & \longrightarrow & W \\ \downarrow & & \downarrow \\ \mathbf{G}^2_{m_k} & \xrightarrow{\varphi} & \mathbf{G}^2_{m_k} \end{array}$$

On en déduit immédiatement l'isomorphisme canonique $\varphi^\star(\mathscr{F} \boxtimes \mathscr{F}) \simeq \mathscr{F} \boxtimes \mathscr{F}$ voulu. Par construction, $\varphi^\star(\mathscr{F} \boxtimes \mathscr{F}) \simeq \mathscr{F}_1 \otimes \mu^\star\mathscr{F}$. On en déduit par récurrence sur $n$ que le faisceau $\mathscr{F}^{\boxtimes n}$ sur $\mathbf{G}_{m_k}^{\,n}$ s'identifie à $\mathscr{F}_1^{\otimes n-1} \otimes m^\star\mathscr{F}$ où $m\colon \mathbf{G}^n_{m_k} \to \mathbf{G}_m$ est le morphisme de multiplication. Ceci achève la démonstration du lemme. $\square$

Montrons la proposition 1.5.2. Soit $n \geq 1$. Les $n$-fonctions symétriques élémentaires définissent un isomorphisme $(\mathbf{A}^1_k)^{\langle n \rangle} \xrightarrow{\sim} \mathbf{A}^n_k$. Cet isomorphisme induit un isomorphisme $\mathbf{G}_{m_k}{}^{\langle n \rangle} \xrightarrow{\sim} \mathbf{A}^{n-1}_k \times \mathbf{G}_{m_k}$, la projection $\pi : \mathbf{G}_{m_k}{}^{\langle n \rangle} \to \mathbf{G}_{m_k}$ étant induite par le produit $\mathbf{G}_{m_k}{}^n \to \mathbf{G}_{m_k}$. D'après le lemme précédent, il vient que quitte à choisir un générateur de $\mathscr{F}_1$, le faisceau $\mathscr{F}^{\langle n \rangle}$ sur $\mathbf{G}_{m_k}{}^{\langle n \rangle}$ s'identifie à $\pi^\star \mathscr{F}$.

Le faisceau $(j_! \mathscr{F})^{\langle n \rangle}$ sur $(\mathbf{P}^1_k)^{\langle n \rangle} \simeq \mathbf{P}^n_k$ s'identifie au prolongement par zéro du faisceau $\mathscr{F}^{\langle n \rangle}$ sur $\mathbf{G}_{m_k}{}^{\langle n \rangle}$. Le complémentaire de $\mathbf{G}_{m_k}{}^{\langle n \rangle}$ dans $(\mathbf{P}^1_k)^{\langle n \rangle}$ s'identifie à la réunion de deux hyperplans $H_0$ et $H_\infty$ qui s'intersectent transversalement. La description donnée plus haut du faisceau localement constant $\mathscr{F}^{\langle n \rangle}$ sur $\mathbf{G}_{m_k}{}^{\langle n \rangle}$ montre qu'il est modéré le long de $H_0$ et $H_\infty$. D'après [Saito 2017b, theorem 7.14], on obtient l'égalité suivante :

$$\begin{aligned} CC(\mathbf{P}^n_k, (j_! \mathscr{F})^{\langle n \rangle}) &= (-1)^n \left( [T^\star_{\mathbf{P}^n_k} \mathbf{P}^n_k] + [T^\star_{H_0} \mathbf{P}^n_k] + [T^\star_{H_\infty} \mathbf{P}^n_k] + [T^\star_{H_0 \cap H_\infty} \mathbf{P}^n_k] \right) \\ &= (-1)^n \left( \tau^\star_{(n)} + \tau^\star_{[0],(n-1)} + \tau^\star_{[\infty],(n-1)} + \tau^\star_{[0]+[\infty],(n-2)} \right) \end{aligned}$$

(Si $n = 1$, $H_0 \cap H_\infty$ est vide, donc le dernier terme de la dernière égalité n'a pas de sens et doit être enlevé.) Ceci montre que l'on a bien l'égalité :

$$S_{j_! \mathscr{F}} = \sum_{n=0}^{\infty} CC(\mathbf{P}^n_k, (j_! \mathscr{F})^{\langle n \rangle}) = \left( \sum_{n=0}^{\infty} (-1)^n \tau^\star_{(n)} \right) \cdot S_{j_! \Lambda} = S_\Lambda \cdot (1 - \tau^\star_{[0]}) \cdot (1 - \tau^\star_{[\infty]}) = S_{j_! \Lambda}.$$

**1.5.4. Proposition.** *Soit $k$ un corps algébriquement clos. Soit $\Lambda$ un corps fini de caractéristique inversible dans $k$. Soit $\mathscr{F}$ un faisceau constructible de $\Lambda$-modules sur $\mathbf{P}^1_k$. On suppose que la restriction de $\mathscr{F}$ à $\mathbf{G}_{m_k}$ est un faisceau localement constant de rang $r$ qui devient trivial après image inverse par le morphisme $[N] : \mathbf{G}_{m_k} \to \mathbf{G}_{m_k}$ d'élévation à la puissance $N$, où $N$ est un entier inversible dans $k$. On suppose aussi que le morphisme canonique $\mathscr{F} \to j_\star j^\star \mathscr{F}$ est injectif, où $j : \mathbf{G}_{m_k} \to \mathbf{P}^1_k$ est l'immersion ouverte canonique. On suppose enfin que $\mathscr{F}_\infty = 0$. Notons $a_0 := r - \dim_\Lambda \mathscr{F}_0$ et $a_\infty = r$. Alors,*

$$S_{\mathscr{F}} = S_{\Lambda^r} \cdot (1 - \tau^\star_{[0]})^{a_0} \cdot (1 - \tau^\star_{[\infty]})^{a_\infty}.$$

*De plus, pour tout entier $n \in \mathbf{N}$, le support singulier de $\mathscr{F}^{\langle n \rangle}$ est exactement le support de son cycle caractéristique.*

*Démonstration.* L'hypothèse fait que $\mathscr{F}$ possède un sous-faisceau $\mathscr{G}$ qui est constant sur $\mathbf{A}^1_k$ de valeur $\mathscr{F}_0$. Notons $\mathscr{H}$ le quotient $\mathscr{H} := \mathscr{F}/\mathscr{G}$. D'après la proposition 1.4.11, la série $S_{\mathscr{G}}$ est connue. De plus, si on admet provisoirement que le résultat voulu est vrai pour $\mathscr{H}$ qui vérifie de plus $\mathscr{H}_0 = 0$, alors on peut appliquer la proposition 1.2.11, et la formule $S_{\mathscr{F}} = S_{\mathscr{G}} \cdot S_{\mathscr{H}}$ permet de conclure.

On peut donc supposer que les deux fibres $\mathscr{F}_0$ et $\mathscr{F}_\infty$ sont nulles. Le $\Lambda$-module $\mathscr{F}_1$ est muni d'une action du groupe abélien $\Gamma = \mu_N(k)$ et l'endomorphisme de $\mathscr{F}_1$ associé à un générateur du groupe $\Gamma$ devient trigonalisable sur une extension finie $\Lambda'$ du corps des coefficients $\Lambda$. Quitte à remplacer $\mathscr{F}$ par $\Lambda' \otimes_\Lambda \mathscr{F}$ et $\Lambda'$ par $\Lambda$, on peut supposer qu'il existe une filtration de $\mathscr{F}_1$ comme $\Lambda[\Gamma]$-module dont les gradués successifs soient de dimension 1. À cette filtration de $\mathscr{F}_1$ correspond une filtration de $\mathscr{F}$ dont les gradués successifs $\mathscr{L}_i$ pour $i \in \{1, \ldots, r\}$ sont de rang 1 sur $\mathbf{G}_{m_k}$. Pour tout $i \in \{1, \ldots, r\}$, on a l'égalité de séries $S_{\mathscr{L}_i} = S_{j_! \Lambda}$ d'après la proposition 1.5.2, on a alors :

$$\begin{aligned} S_{\mathscr{F}} &= \prod_{i=1}^{r} S_{\mathscr{L}_i} = (S_{j_! \Lambda})^r = S_{j_! \Lambda^r} \\ &= S_{\Lambda^r} \cdot (1 - \tau^\star_{[0]})^r \cdot (1 - \tau^\star_{[\infty]})^r. \end{aligned}$$

□

Nous pouvons maintenant démontrer le théorème 1.5.1. On procède comme en 1.3.12 et comme dans la démonstration de 1.4.8. Soient $n \in \mathbf{N}$ et $D \in X^{\langle n \rangle}(k)$. Il s'agit de montrer que la formule énoncée dans le théorème pour le cycle caractéristique et le support singulier de $\mathscr{F}^{\langle n \rangle}$ est vraie au voisinage de $D$. On écrit $D = \sum_{j \in J} \nu_j x_j$ avec des points $x_j \in X(k)$ distincts. Soit $j \in J$. L'hypothèse faite sur $\mathscr{F}$ au voisinage de $x_j$ permet de choisir $f_j : U_j \to X$ étale telle que $f_j^{-1}(x_j) = \{u_j\}$, $g_j : U_j \to \mathbf{P}^1_k$ étale telle que $g_j^{-1}(0) = \{u_j\}$ et $g_j^{-1}(\infty) = \varnothing$ et $\mathscr{G}_j$ un faisceau sur $\mathbf{P}^1_k$ tel que $f_j^\star \mathscr{F} \simeq g_j^\star \mathscr{G}_j$, le faisceau $\mathscr{G}_j$ vérifiant les hypothèses de la proposition 1.5.4.

Notons $U := \coprod_{j \in J} U_j$ et $Y := \coprod_{j \in J} \mathbf{P}^1_k$. De même que pour le lemme 1.3.11, en utilisant les morphismes étales évidents $U \to X$ et $U \to Y$, pour montrer la formule du théorème pour $\mathscr{F}^{\langle n \rangle}$ au voisinage de $D$, il suffit de la montrer pour $\mathscr{F}_{|U}{}^{\langle n \rangle}$ au voisinage de $D_U := \sum_{j \in J} \nu_j u_j$, et pour cela il suffit de la montrer au voisinage de l'image $D_Y := \sum_{j \in J} \nu_j y_j$ du cycle $D_U$ par le morphisme évident $U \to Y$, et ce pour le faisceau $\mathscr{G}^{\langle n \rangle}$ où $\mathscr{G}$ est un faisceau sur $Y$ ayant la propriété d'induire $\mathscr{G}_j$ sur l'ouvert $Y_j \simeq \mathbf{P}^1_k$ de $Y$ pour tout $j \in J$.

Grâce à la proposition 1.5.4, on peut calculer le cycle caractéristique et le support singulier de $\mathscr{G}_j^{\langle m \rangle}$ sur $Y_j^{\langle m \rangle}$ pour tout $m \in \mathbf{N}$ :

$$S_{\mathscr{G}_j} = \left(S_\Lambda\right)^r \cdot (1 - \tau_{y_j}^\star)^{a_j}$$

où on a noté $a_j$ la chute du rang en $y_j$. Cette identité correspond à des égalités de cycles dans $T^\star Y_j^{\langle m \rangle}$ pour tout $m \in \mathbf{N}$. En identifiant, $\mathscr{G}_j$ à un faisceau sur $Y$ qui s'annule sur tous les ouverts $Y_{j'}$ pour $j' \neq j$, on peut récrire cette identité en termes de cycles dans $T^\star Y^{\langle m \rangle}$ pour tout $m \in \mathbf{N}$ :

$$S_{\mathscr{G}_j} = \left(S_{\Lambda_{Y_j}}\right)^r \cdot (1 - \tau_{y_j}^\star)^{a_j} \cdot (1 - \tau_{z_j}^\star)^r$$

où pour tout $j \in J$, $y_j$ et $z_j$ sont les deux points de $Y_j \simeq \mathbf{P}^1_k$ correspondant respectivement à 0 et $\infty$. Avec les identifications précédentes, on a $\mathscr{G} \simeq \oplus_{j \in J} \mathscr{G}_j$, donc en appliquant la proposition 1.2.10, on obtient :

$$\begin{aligned} S_{\mathscr{G}} &= \prod_{j \in J} S_{\mathscr{G}_j} \\ &= \left(\prod_{j \in J} S_{\Lambda_{Y_j}}\right)^r \cdot \prod_{j \in J} (1 - \tau_{y_j}^\star)^{a_j} \cdot \prod_{j \in J} (1 - \tau_{z_j}^\star)^r \\ &= (S_\Lambda)^r \cdot \prod_{j \in J} (1 - \tau_{y_j}^\star)^{a_j} \cdot \prod_{j \in J} (1 - \tau_{z_j}^\star)^r \\ &= S_{\Lambda^r} \cdot \prod_{j \in J} (1 - \tau_{y_j}^\star)^{a_j} \cdot \prod_{j \in J} (1 - \tau_{z_j}^\star)^r \end{aligned}$$

Pour l'avant-dernière égalité, on a utilisé la décomposition $\Lambda \simeq \oplus_{j \in J} \Lambda_{Y_j}$. Notons $Y' \subset Y$ l'ouvert $\coprod_{j \in J} \mathbf{A}^1_k \subset \coprod_{j \in J} \mathbf{P}^1_k$. En restreignant l'identité précédente aux ouverts $T^\star Y'^{\langle m \rangle}$ de $T^\star Y^{\langle m \rangle}$ pour tout $m \in \mathbf{N}$, on obtient finalement :

$$S_{\mathscr{G}_{|Y'}} = S_{\Lambda^r} \cdot \prod_{j \in J} (1 - \tau_{y_j}^\star)^{a_j}$$

En passant à la composante de degré $n$, on obtient le résultat voulu au voisinage de $D_Y \in T^\star Y^{\langle n \rangle}$, ce qui achève la démonstration du théorème.

## 2. Morphisme d'Abel-Jacobi et théorèmes d'acyclicité de Deligne

L'objectif de cette section est d'établir un théorème d'acyclicité (**2.1**), dû à Pierre Deligne, en utilisant le théorème principal de la section précédente. Plus précisément, on va montrer que pour $\mathscr{F}$ comme dans le théorème **1.5.1**, la courbe $X$ *propre*, et l'entier $n$ suffisamment grand (de façon explicite), la paire $(\int^n : X^{\langle n\rangle} \to \mathrm{Pic}^n_X, \mathscr{F}^{\langle n\rangle})$, où $\int^n$ est le *morphisme d'Abel-Jacobi* (annexe B), est localement acyclique. De plus, pour $n$ « critique », cette paire est localement acyclique au-dessus du complémentaire d'un point explicite de $\mathrm{Pic}^n_X$. Dans la section suivante, nous en déduisons, toujours selon Pierre Deligne, que le déterminant de la cohomologie de $X$ à valeurs dans $\mathscr{F}$ se « localise ».

**2.1. Théorème** (Pierre Deligne, *circa* 1980)**.** *Soient $X$ une courbe projective lisse sur un corps algébriquement clos $k$, $\Lambda$ un anneau fini de torsion inversible dans $k$ et $\mathscr{F}$ un $\Lambda$-faisceau constructible sur $X$. On suppose qu'il existe une immersion ouverte dominante $j : U \hookrightarrow X$ et un entier $r$ tel que $j^\star\mathscr{F}$ soit localement constant, localement libre de rang $r$, modéré le long de $S := X - U$, localement libre sur $S$, et que le morphisme d'adjonction $\mathscr{F} \to j_\star j^\star\mathscr{F}$ soit injectif. On suppose de plus que pour tout $s \in \mathscr{F}$ la fibre $\mathscr{F}_s$ est un $\Lambda$-module libre. On note alors $\Delta_{\mathscr{F}} := \sum_{s\in S} a_s \cdot s$, où $a_s := r - \dim_\Lambda \mathscr{F}_s$ est la chute du rang en $s$, puis $n_{\mathscr{F}} := r(2g-2) + \deg(\Delta_{\mathscr{F}})$ et, lorsque $n_{\mathscr{F}} > 0$, posons enfin $K_{\mathscr{F}} := [\Delta_{\mathscr{F}}] + rK_X \in \mathrm{Pic}^{n_{\mathscr{F}}}_X$. Alors, la paire $(\int^n : X^{\langle n\rangle} \to \mathrm{Pic}^n_X, \mathscr{F}^{\langle n\rangle})$ est* localement acyclique *pour tout entier* naturel *$n > n_{\mathscr{F}}$ et, si $n = n_{\mathscr{F}}$, localement acyclique au-dessus de $\mathrm{Pic}^{n\times}_X := \mathrm{Pic}^n_X - \{K_{\mathscr{F}}\}$.*

La méthode que allons suivre est celle de Deligne : la seule différence est dans notre usage du support singulier. Rappelons ([Beilinson 2017, §1.3]) qu'une paire $(f : Y \to Z, \mathscr{G})$ est localement acyclique si la préimage par $df^\vee : T^\star Z_Y := T^\star Z \times_Z Y \to T^\star Y$ du fermé $SS(\mathscr{G})$ est contenue dans la section nulle.

**2.1.1.** Soient $X$, $\mathscr{F}$ et $n$ comme en **1.5.1**. Voyons à quelles conditions la paire $(\mathscr{F}^{\langle n\rangle}, \int^n : X^{\langle n\rangle} \to \mathrm{Pic}^n_X)$ est acyclique, sous l'hypothèse supplémentaire que $X$ est *propre*. D'après ce qui précède (appliqué fibre à fibre), il suffit que pour $\Delta, \underline{e}$ tels que $\tau^\star_{\Delta,e}$ apparaisse avec un coefficient non nul dans l'expression de $CC(X^{\langle n\rangle}, \mathscr{F}^{\langle n\rangle})$, on ait, pour tout diviseur effectif $D$ de degré $n$ correspondant à un point géométrique de $X^{\langle n\rangle}$ :

$$(\dagger) \qquad \left\{\omega \in (T^\star \mathrm{Pic}^n_X)_{[D]} : d\textstyle\int^{n\vee}_D(\omega) \in (\tau^\star_{\Delta,\underline{e}})_D \subseteq (T^\star X^{\langle n\rangle})_D\right\} = \{0\}.$$

D'après les résultats classiques rappelés en **B.3**, l'application $d\int^{n\vee}_D$ entre les fibrés cotangents s'identifie au morphisme

$$H^0(X, \Omega^1_X) \to H^0(X, \mathscr{O}_X(D)/\mathscr{O}_X)^\vee$$
$$\omega \mapsto \mathrm{rés}(\text{—} \cdot \omega).$$

Fixons $\underline{e}$ comme en **1.5.1**, en supposant tout d'abord $\Delta = 0$ pour simplifier. Rappelons (**1.1.2**) que si $D$ se décompose en $D_1 + 2D_2 + \cdots + nD_n$, la fibre en $D$ du fibré $\tau^\star_{\underline{e}}$ est, par définition, l'orthogonal dans $H^0(\mathscr{O}_X(D_1 + 2D_2 + \cdots + nD_n)/\mathscr{O}_X)^\vee$ de $H^0(\mathscr{O}_X(D_1 + D_2 + \cdots + D_n)/\mathscr{O}_X) \subseteq H^0(\mathscr{O}_X(D_1 + 2D_2 + \cdots + nD_n)/\mathscr{O}_X)$, c'est-à-dire l'espace des formes linéaires nulles sur $\mathscr{O}_{D^{\downarrow\underline{e}}}(D^{\downarrow\underline{e}}) \subseteq \mathscr{O}_D(D)$, où on pose $D^{\downarrow\underline{e}} := D_1 + D_2 + \cdots + D_n$. Ainsi, le terme de gauche de $(\dagger)$ (pour $\Delta = 0$) est en bijection avec l'ensemble des formes différentielles $\omega$ telles que $\mathrm{rés}(\mathscr{O}_{D^{\downarrow\underline{e}}}(D^{\downarrow\underline{e}}) \cdot \omega) = 0$, autrement dit $\omega \in H^0(X, \Omega^1_X(-D^{\downarrow\underline{e}}))$. Il s'agit de montrer qu'un tel espace est nul.
Plus précisément, prenant maintenant en compte la translation par un diviseur effectif $\Delta \leqslant \Delta_{\mathscr{F}}$, on a acyclicité de la paire considérée si :

> pour tous $\Delta, \underline{e}$ tels que $\tau^\star_{\Delta,\underline{e}}$ apparaisse avec un coefficient non nul dans l'expression de $CC(X^{\langle n\rangle}, \mathscr{F}^{\langle n\rangle})$ en **1.5.1**, et toute décomposition de $D \in X^{\langle n\rangle}$ en somme $\Delta + D'$ de diviseurs effectifs, on a

$$H^0(X, \Omega^1_X(-D'^{\downarrow \underline{e}})) = 0.$$

**2.1.2**. Observons (Riemann-Roch) que si $E$ est un diviseur effectif de degré $> 2g-2$, alors $H^0(X, \Omega^1_X(-E)) = 0$. La paire $(\mathscr{F}^{\langle n\rangle}, \int^n)$ est donc acyclique dès lors que les $D'^{\downarrow \underline{e}}$ ci-dessus sont *tous* de degré $> 2g-2$.

Supposons $n > r(2g-2) + \deg(\Delta_{\mathscr{F}})$. Dans ce cas, tout diviseur $D'$ tel que $D = \Delta + D'$ satisfait $\deg(D') > r(2g-2)$, car $\deg(\Delta) \leqslant \deg(\Delta_{\mathscr{F}})$. D'autre part, on a

$$\deg(D') \leqslant r \deg(D'^{\downarrow \underline{e}})$$

car les coefficients binomiaux $\binom{r}{i}$ sont nuls pour $i > r$ : le diviseur $D'$ est une somme $D'_1 + 2D'_2 + \cdots + rD'_r$ (avec $D'_i$ de degré $i \leqslant r$). Ces deux inégalités entraînent $\deg(D'^{\downarrow \underline{e}}) > 2g-2$.

Supposons $n = r(2g-2) + \deg(\Delta_{\mathscr{F}})$. Dans ce cas, une condition nécessaire pour qu'il n'y ait pas acyclicité au-dessus de $D$ est que $\deg(\Delta) = \deg(\Delta_{\mathscr{F}})$ — donc $\Delta = \Delta_{\mathscr{F}}$, car $\Delta \leqslant \Delta_{\mathscr{F}}$ —, $D' = rD'^{\downarrow \underline{e}}$ et que $[D'^{\downarrow \underline{e}}]$ soit la classe $K_X$ du diviseur canonique. Autrement dit : $D = \Delta_{\mathscr{F}} + rD''$, où $[D''] = K_X$. (On a alors $e_r = 2g-2$, les autres entiers étant nuls.)

**2.1.3**. Ceci achève la démonstration du théorème, du moins dans le cas particulier où l'anneau $\Lambda$ de coefficients est un corps.

*Stricto sensu*, le théorème précédent n'a été établi que pour $\Lambda$ un corps fini de caractéristique inversible sur $k$. Voyons rapidement comment le cas général s'en déduit. Tout d'abord, on peut supposer $\Lambda$ connexe. Son réduit $\overline{\Lambda}$ est donc, sous notre hypothèse, un corps fini. Notons $\mathscr{G}$ le faisceau $\mathscr{F}^{\langle n\rangle}$ et $\overline{\mathscr{G}}$ le produit tensoriel $\mathscr{G} \otimes_\Lambda \overline{\Lambda}$ qui est isomorphe à $\left(\mathscr{F} \otimes_\Lambda \overline{\Lambda}\right)^{\langle n\rangle}$ d'après [Roby 1963, théorème III.3]. Pour établir l'acyclicité de $(\int^n, \mathscr{G})$, il suffit de vérifier celle de $(\int^n, \overline{\mathscr{G}})$ car d'une part $\mathscr{G}$ est muni d'une filtration finie dont les gradués sont isomorphes à $\overline{\mathscr{G}}$ et, d'autre part, l'acyclicité est stable par extension de faisceaux. Le rang et les chutes du rang aux points de $S$ étant inchangés, il suffit pour conclure de vérifier que le faisceau de $\overline{\Lambda}$-modules $\overline{\mathscr{G}}$ satisfait également les hypothèses du théorème : modération et injectivité du morphisme de spécialisation. C'est évident pour la première : si un $P$-Sylow agit trivialement, il en est de même sur tout quotient. Quant à la seconde, elle se ramène par passage aux fibres en les points de $S$ à montrer que si $F \to L$ est un morphisme *injectif* de $\Lambda$-modules libres de type fini, il en est de même de $\overline{F} \to \overline{L}$. Cela résulte du fait que, l'anneau $\Lambda$ étant local artinien, $F$ est nécessairement un facteur direct de $L$. Ceci achève la démonstration du théorème avec le niveau de généralité annoncé.

**2.2**. **Projection sur la droite affine.** Soient $X$ et $\mathscr{F}$ comme dans l'énoncé du théorème 2.1, dont on reprend les notations. Soit $t : (\mathsf{Pic}^n_X, K) \dashrightarrow (\mathbf{A}^1, 0)$ un germe de morphisme, où $K$ est une classe de diviseur de degré $n$. Notons $f_t : X^{\langle n\rangle} \dashrightarrow \mathbf{A}^1$ le germe de morphisme composé $t \circ \int^n$, défini sur un voisinage ouvert de la fibre $X^{\langle n\rangle}_K$. Ici encore, l'acyclicité de paire $(f_t, \mathscr{F}^{\langle n\rangle})$ en un diviseur effectif $D$ de degré $n$ est garantie dès lors que les préimages analogues à (†) (pour $\Delta$ et $e$ variables) sont toutes nulles, où on précompose $d\int^{n\vee}_D : (T^\star \mathsf{Pic}^n_X)_{[D]} \to (T^\star X^{\langle n\rangle})_D$ par $dt^\vee_D : (T^\star \mathbf{A}^1)_0 \to (T^\star \mathsf{Pic}^n_X)_{[D]}$. Au niveau des groupes de cohomologie cohérente, cela revient à considérer le morphisme composé

$$\begin{aligned} k \to H^0(X, \Omega^1_X) &\to H^0(X, \mathscr{O}_X(D)/\mathscr{O}_X)^\vee \\ 1 \mapsto \omega_t := dt &\mapsto \text{rés}(\text{—} \cdot \omega_t), \end{aligned}$$

où $dt$ est la différentielle de $t$ en $K$.

Comme ci-dessus, il ne peut donc y avoir de singularité de la paire en $D$ que s'il existe une décomposition en diviseurs effectifs $D = \Delta + D'$ avec $\Delta \leqslant \Delta_{\mathscr{F}}$ telle que

$$0 \neq \omega_t \in H^0(X, \Omega^1_X(-D'^{\downarrow \underline{e}})),$$

c'est-à-dire $\mathrm{div}(\omega_t) \geqslant D'^{\downarrow \underline{e}}$. On l'a vu, le groupe des sections globales de droite est nul si $n > n_{\mathscr{F}}$, d'où l'acyclicité en tout point de $X^{\langle n \rangle}_K$ dans ce cas. Si $n = n_{\mathscr{F}}$, il est non nul si et seulement si $D = \Delta_{\mathscr{F}} + rD''$ avec $[D''] = K_X$ (de sorte que $K = K_{\mathscr{F}}$) et la condition $\omega_t \in H^0(X, \Omega^1_X(-D''))$ se réécrit $D'' = \mathrm{div}(\omega_t)$.

**2.3. Théorème** (Pierre Deligne). *Soient $X$ et $\mathscr{F}$ comme dans l'énoncé du théorème 2.1, un entier $n \geqslant n_{\mathscr{F}}$, une fonction $t : \mathrm{Pic}^n_X \dashrightarrow \mathbf{A}^1_k$ définie sur un ouvert $W$ de $\mathrm{Pic}^n_X$ de différentielle en $K_{\mathscr{F}}$ notée $\omega_t := dt$, vue comme forme différentielle sur $X$. Alors, la paire $\left(t \circ \int^n : X^{\langle n \rangle}_W \to \mathbf{A}^1_k, \mathscr{F}^{\langle n \rangle}\right)$ est* localement acyclique *sauf peut-être, lorsque $n = n_{\mathscr{F}}$, en le point $\Delta_{\mathscr{F},\omega_t} := \Delta_{\mathscr{F}} + r\mathrm{div}(\omega_t)$ au-dessus de $K_{\mathscr{F}} \in \mathrm{Pic}^{n_{\mathscr{F}}}_X$.*

## 3. Localité du facteur epsilon

**3.1. Lemme.** *Soient $A$ un torseur sous une variété abélienne de dimension $\geqslant 2$ sur un corps algébriquement clos $k$, un point rationnel $a \in A(k)$, un anneau $\Lambda$ de coefficients fini de cardinal inversible dans $k$, et $\mathscr{H} \in \mathrm{D}^{\mathrm{b}}_{\mathrm{ctf}}(A, \Lambda)$ un complexe dont la restriction à l'ouvert $A^\times := A - \{a\}$ est à objets de cohomologie localement constants, à fibres projectives de type fini. Pour tout germe de fonction $t : (A, a) \dashrightarrow (\mathbf{A}^1_k, 0)$ lisse en $a$, il existe un isomorphisme canonique*

$$\mathrm{dét}_\Lambda(\mathrm{R}\Gamma(A, \mathscr{H})) \simeq \mathrm{dét}_\Lambda^{-1} \Phi_t(\mathscr{H})_a \, .$$

D'après les idées introduites dans [Deligne 1987], on utilise ici la notation $\mathrm{dét}_\Lambda K$ pour le déterminant d'un complexe parfait $K$ de $\Lambda$-modules. On utilisera le fait que si $K' \to K \to K'' \to K'[1]$ est un triangle distingué, on dispose d'un isomorphisme *canonique* $\mathrm{dét}_\Lambda K \simeq \mathrm{dét}_\Lambda K' \otimes \mathrm{dét}_\Lambda K''$, cf. [Knudsen et Mumford 1976] et [Muro, Tonks et Witte 2015] pour plus de détails.

*Démonstration.* Commençons par traiter le cas particulier où $\mathscr{H}$ est un *faisceau*, c'est-à-dire est concentré en un seul degré, que l'on peut supposer être 0. Soit

$$\mathscr{H} \to R^0 j_\star j^\star \mathscr{H} =: \overline{\mathscr{H}} \to \mathscr{Q}$$

le triangle distingué canonique, où $j$ est l'immersion ouverte $A^\times \hookrightarrow A$ et la première flèche est l'unité de l'adjonction. L'hypothèse que $a$ est de codimension au moins 2 dans $A$ entraîne, par pureté, que le faisceau $\overline{\mathscr{H}}$ est *lisse* sur $A$. D'autre part, le *complexe* $\mathscr{Q}$ est à support dans $\{a\}$. D'après [Deligne 1974, construction 1, p. 70], on dispose d'un isomorphisme canonique $\mathrm{dét}_\Lambda(\mathrm{R}\Gamma(A, \overline{\mathscr{H}})) \xrightarrow{\sim} \Lambda$ de sorte que le triangle ci-dessus induit un isomorphisme canonique $\mathrm{dét}_\Lambda(\mathrm{R}\Gamma(A, \mathscr{H})) \xrightarrow{\sim} \mathrm{dét}_\Lambda^{-1} \mathscr{Q}_a$. (Rappelons que, sous nos hypothèses, le complexe $\mathrm{R}\Gamma(A, \mathscr{H})$ est *parfait.*) D'autre part, notant $V$ l'ouvert contenant $a$ sur lequel le morphisme $t$ est défini, et $i$ l'immersion fermée $V_0 \hookrightarrow V$, la nullité des cycles proches $\Psi_t \mathscr{Q}_{|V}$ du complexe à support ponctuel $\mathscr{Q}$ induit un isomorphisme $\Phi_t(\mathscr{Q}_{|V}) \xrightarrow{\sim} i^\star \mathscr{Q}_{|V}[1]$. Enfin, le triangle distingué $\overline{\mathscr{H}} \to \mathscr{Q} \to \mathscr{H}[1]$ induit par acyclicité de $(t, \overline{\mathscr{H}}_{|V})$ un isomorphisme $\Phi_t(\mathscr{Q}_{|V}) \xrightarrow{\sim} \Phi_t \mathscr{H}_{|V}[1]$ d'où un isomorphisme $i^\star \mathscr{Q}_{|V} \xrightarrow{\sim} \Phi_t \mathscr{H}_{|V}$ et le résultat annoncé en passant au déterminant des fibres en $a$.

Cas d'un complexe. Pour tout entier relatif $n$, le triangle distingué

$$H^n(\mathscr{H})[-n] \to \tau_{\geqslant n} \mathscr{H} \to \tau_{>n} \mathscr{H}$$

induit par application des foncteurs triangulés $F_0 := \mathrm{R}\Gamma(A, -)$ et $F_1 := \Phi_t(-_{|V})_a$ deux triangles distingués, pour $i \in \{0, 1\}$ :

$$F_i(H^n(\mathscr{H})[-n]) \to F_i(\tau_{\geqslant n} \mathscr{H}) \to F_i(\tau_{>n} \mathscr{H}).$$

Passant au déterminant de ces complexes parfaits, on obtient pour chaque $n$ deux isomorphismes canoniques $\text{dét}_\Lambda F_i(\tau_{\geqslant n}\mathscr{H}) \xrightarrow{\sim} \text{dét}_\Lambda F_i(H^n\mathscr{H})\otimes_\Lambda \text{dét}_\Lambda F_i(\tau_{>n}\mathscr{H})$ d'où, par récurrence sur l'amplitude de $\mathscr{H}$, un diagramme d'isomorphismes

$$\begin{array}{ccc} \text{dét}_\Lambda \mathrm{R}\Gamma(A,\tau_{\geqslant n}\mathscr{H}) & \longrightarrow & \text{dét}_\Lambda \mathrm{R}\Gamma(A,H^n\mathscr{H}) \otimes_\Lambda \text{dét}_\Lambda \mathrm{R}\Gamma(A,\tau_{>n}\mathscr{H}) \\ \downarrow & & \downarrow \qquad\qquad \downarrow \\ \text{dét}_\Lambda \Phi_t(\tau_{\geqslant n}\mathscr{H}_{|V})_a & \longrightarrow & \text{dét}_\Lambda \Phi_t(H^n\mathscr{H}_{|V})_a \otimes_\Lambda \text{dét}_\Lambda \Phi_t(\tau_{>n}\mathscr{H}_{|V})_a. \end{array}$$

Pour $n$ suffisamment grand, on obtient l'énoncé désiré. □

**3.2. Lemme.** *Soient $S$ un trait strictement hensélien, $(Y,y) \xrightarrow{f} (A,a) \overset{t}{\dashrightarrow} (S,s)$ des morphismes de type fini pointés en des points fermés, et $\mathscr{G}$ un faisceau abélien de torsion sur $Y$. Soit $V$ l'ouvert de $A$ sur lequel $t$ est défini et $f_t$ le morphisme composé $Y_V \to V \to S$. Si la paire $(f_t,\mathscr{G})$ est* localement acyclique *hors de $\{y\}$ et $f$ est* propre, *le morphisme canonique*

$$\Phi_t(\mathrm{R}f_\star\mathscr{G})_a \to \Phi_{f_t}(\mathscr{G})_y$$

*est un* isomorphisme.

*Démonstration.* On peut supposer $V = A$. Par propreté de $f$, le morphisme $\Phi_t(\mathrm{R}f_\star\mathscr{G}) \to \mathrm{R}f_{s\star}\Phi_{f_t}(\mathscr{G})$ déduit du morphisme analogue pour les cycles proches ([SGA 7 XIII, 2.1.7.1]) est un isomorphisme ; passant à la fibre en $a$, on obtient $\Phi_t(\mathrm{R}f_\star\mathscr{G})_a \xrightarrow{\sim} R\Gamma(Y_a,\Phi_{f_t}(\mathscr{G}))$. La conclusion résulte alors du fait que $\Phi_{f_t}(\mathscr{G})$ est supporté en $y \in Y_a$. □

**3.3. Théorème** (Pierre Deligne, *circa* 1980)**.** *Soient $X$ une courbe projective lisse sur un corps algébriquement clos $k$ de caractéristique différente d'un nombre premier $\ell$, un corps de coefficients $\Lambda$ égal à $\mathbb{F}_\ell$ ou $\mathbb{Q}_\ell$, un $\Lambda$-faisceau constructible $\mathscr{F}$ sur $X$ satisfaisant les hypothèses du théorème 2.1 (dont la modération), $n := \chi(X,\mathscr{F}[1]) = -\chi(X,\mathscr{F})$, supposé strictement positif, sa caractéristique d'Euler-Poincaré (au signe près), et $\Delta_\mathscr{F}$ (resp. $K_\mathscr{F}$) le diviseur effectif associé (resp. la classe $rK_X + [\Delta_\mathscr{F}] \in \mathrm{Pic}^n_X$). Alors, pour tout germe de fonction lisse $t : (\mathrm{Pic}^n_X, K_\mathscr{F}) \dashrightarrow (\mathbf{A}^1_k, 0)$, de différentielle en $K_\mathscr{F}$ la forme notée $\omega_t \in H^0(X,\Omega^1_{X/k})$, on a un isomorphisme* canonique

$$\varepsilon(X,\mathscr{F}) := \text{dét}_\Lambda \mathrm{R}\Gamma(X,\mathscr{F}) \xrightarrow{\sim} \text{dét}_\Lambda^{(-1)^n} \Phi_{\int_t^n}\big(\mathscr{F}^{\langle n\rangle}\big)_{\Delta_{\mathscr{F},\omega}},$$

*où $\Delta_{\mathscr{F},\omega_t} \in X^{\langle n\rangle}$ est le diviseur effectif $\Delta_\mathscr{F} + r\mathrm{div}(\omega)$ et $\int_t^n$ est le composé du morphisme d'Abel-Jacobi $(X^{\langle n\rangle},\Delta_{\mathscr{F},\omega}) \to (\mathrm{Pic}^n_X, K_\mathscr{F})$ avec le germe de fonction $t$.*

*Démonstration.* Soit $K$ le complexe parfait $\mathrm{R}\Gamma(X,\mathscr{F}[1])$. Trivialement, on a $\varepsilon(X,\mathscr{F}) \xrightarrow{\sim} \text{dét}_\Lambda^{-1} K$. Comme d'autre part $K$ est de caractéristique d'Euler-Poincaré $n > 0$, on a d'après [Deligne 1974, b.4, cor. 1] $\text{dét}\, K \xrightarrow{\sim} \text{dét}\,(\mathrm{L}\Lambda^n K)$, où $\mathrm{L}\Lambda^n$ est le dérivé au sens de [Dold et Puppe 1961] du foncteur non additif $\Lambda^n$. D'après la formule du décalage de Quillen-Illusie ([Quillen 1968, prop. 7.21], [Illusie 1971-1972, I.4.3.2.1] — déjà utilisée page 20), on a $\mathrm{L}\Lambda^n K \xrightarrow{\sim} \mathrm{L}\Gamma^n\big(\mathrm{R}\Gamma(X,\mathscr{F})\big)[n]$, où $\Gamma^n$ est le foncteur des puissances divisées (cf. rappels en A.1.1). Enfin, la formule de Künneth symétrique ([SGA 4 XVII, 5.5.21]) produit un isomorphisme $\mathrm{L}\Gamma^n\big(\mathrm{R}\Gamma(X,\mathscr{F})\big) \xrightarrow{\sim} \mathrm{R}\Gamma(X^{\langle n\rangle},\mathscr{F}^{\langle n\rangle})$. Par composition, on a construit un isomorphisme

$$\varepsilon(X,\mathscr{F}) \xrightarrow{\sim} \text{dét}_\Lambda\big(\mathrm{R}\Gamma(X^{\langle n\rangle},\mathscr{F}^{\langle n\rangle})[n-1]\big).$$

Les deux lemmes précédents, appliqués à $(f : Y \to A) := (\int^n : X^{\langle n\rangle} \to \mathrm{Pic}^n_X)$, $\mathscr{G} := \mathscr{F}^{\langle n\rangle}$ et $\mathscr{H} := \mathrm{R}\int^n_\star \mathscr{G}$, en les points $a := K_\mathscr{F}$ et $y := \Delta_{\mathscr{F},\omega}$ permettent de réécrire successivement $\text{dét}_\Lambda \mathrm{R}\Gamma(X^{\langle n\rangle},\mathscr{F}^{\langle n\rangle}) = \text{dét}_\Lambda \mathrm{R}\Gamma(\mathrm{Pic}^n_X,\mathscr{H})$ en $\text{dét}_\Lambda^{-1}\Phi_t(\mathscr{H})_{K_\mathscr{F}}$ puis $\Phi_t(\mathscr{H})_{K_\mathscr{F}}$ en une fibre, *locale en* *haut*, de cycles évanescents, comme dans l'énoncé du théorème. Les hypothèses d'acyclicités sont satisfaites en vertu de 2.1 et 2.3. □

**3.4. Corollaire.** *Fixons un entier r, un ouvert dense U de la courbe X, un diviseur effectif Δ à support dans* $S := X - U$, *et* $n := r(2g-2) + \deg(\Delta) > 0$. *Fixons également un germe de fonction lisse* $t : (\mathrm{Pic}^n_X, K_{r,\Delta} := rK_X + [\Delta]) \dashrightarrow (\mathbf{A}^1_k, 0)$ *et posons* $\Sigma := S \cup \mathrm{div}(\omega_t)$, *où* $\omega_t$ *est la différentielle de t en* $K_{r,\Delta}$. *Soient* $\mathscr{F}_1$ *et* $\mathscr{F}_2$ *deux faisceaux sur X comme ci-dessus, ayant les mêmes invariants : même rang r, même ouvert de lissité U, même diviseur* $\Delta_{\mathscr{F}_1} = \Delta = \Delta_{\mathscr{F}_2}$. *Toute famille d'isomorphismes* $\mathscr{F}_{1|X_{(s)}} \xrightarrow{\sim} \mathscr{F}_{2|X_{(s)}}$ *entre les restrictions aux hensélisés en* $s \in \Sigma$ *induit un isomorphisme canonique*

$$\varepsilon(X, \mathscr{F}_1) \xrightarrow{\sim} \varepsilon(X, \mathscr{F}_2),$$

*fonctoriel en un sens évident.*

En effet, d'après le théorème précédent, $\varepsilon(X, \mathscr{F})$ ne dépend que de la restriction de $\mathscr{F}^{\langle n\rangle}$ au hensélisé du produit symétrique en le point $K_{r,\Delta}$.

**3.5. Corollaire** (Pierre Deligne, *circa* 1980)**.** *Soient X et* $\mathscr{F}$ *comme dans l'énoncé du théorème. Pour tout faisceau* $\mathscr{L}$ *lisse de rang 1 sur X, on a un isomorphisme canonique*

$$\varepsilon(X, \mathscr{F} \otimes \mathscr{L}) \simeq \varepsilon(X, \mathscr{F}) \otimes H^0(|K_{\mathscr{F}}|, \mathscr{L}^{\langle n\rangle}),$$

*où* $n := \chi(\mathscr{F}[1]) > 0$ *et* $|K_{\mathscr{F}}| = \int^{n-1}(K_{\mathscr{F}}) \subseteq X^{\langle n\rangle}$ *est le système linéaire associé à* $K_{\mathscr{F}}$.

Voir [Umezaki, Yang et Zhao 2020] pour une généralisation.

## Annexe A. Puissances symétriques : rappels

### A.1. Définition.

**A.1.1.** Soient $A$ un anneau commutatif et $M$ un $A$-module. Pour tout entier $n \geqslant 0$, rappelons que l'on note $\mathrm{TS}^n_A(M)$ le sous $A$-module $\mathrm{Fix}(\mathfrak{S}_n \circlearrowright M^{\otimes_A n})$ de $M^{\otimes_A n} := M \otimes_A \cdots \otimes_A M$ constitué des **tenseurs symétriques de degré** $n$ ([Bourbaki A, chap. IV, §5 n°2], [SGA 4 XVII, 5.5.2.1]) et définit le $A$-module des **puissances divisées** $\mathsf{Y}^n_A(M)$ [②] de degré $n$, représentant le foncteur des applications $n$-iques (au sens de [Roby 1963, I.§8], repris en [SGA 4 XVII, 5.5.2.2]). Ce dernier foncteur satisfait de meilleurs propriétés, notamment la commutation aux changements de base, dont hérite le premier quand on observe que le morphisme naturel $\mathsf{Y}^n_A(M) \to \mathrm{TS}^n_A(M)$, déduit de $\gamma^n : M \to \mathrm{TS}^n_A(M)$, $m \mapsto m \otimes \cdots \otimes m$ par la propriété universelle, est un isomorphisme pour $M$ plat sur $A$. (Il suffit d'ailleurs d'établir ce résultat pour $M$ libre, par passage à la colimite.) Dans le cas général, une présentation plate de $M$ induit une description de $\mathsf{Y}^n_A(M)$ en termes de tenseurs symétriques ([SGA 4 XVII, 5.5.2.5.1]). Le morphisme diagonal induit un morphisme $\mathsf{Y}_A(M) \to \mathsf{Y}_A(M) \otimes_A \mathsf{Y}_A(M)$ ([Roby 1963, V.§2]) ; lorsque $M$ est *plat*, et donc $\mathrm{TS}^r_A(M)$ pour tout $r$ ([SGA 4 XVII, 5.5.2.4]), on en déduit pour toute décomposition $n = n_1 + n_2$ un morphisme $\mathrm{TS}^n_A(M) \to \mathrm{TS}^{n_1}_A(M) \otimes_A \mathrm{TS}^{n_2}_A(M)$, injectif — car le terme de droite n'est alors que le sous-module de $M^{\otimes n}$ des invariants sous $\mathfrak{S}_{n_1} \times \mathfrak{S}_{n_2} \leqslant \mathfrak{S}_n$ (cf. [SGA 4 XVII, 5.5.7.3, démonstration]). Cette structure de cogèbre est définie dans [Bourbaki A, chap. IV, §5 n°7], si $M$ est libre sur $A$.

**A.1.2.** Fixons $n \geqslant 0$, un anneau commutatif $A$ et posons $S := \mathrm{Spec}(A)$. Si $B$ est une $A$-algèbre commutative, il en est de même de $\mathrm{TS}^n_A(B)$ [③] qui représente le $S$-schéma quotient $(\mathrm{Spec}(B)/S)^{\langle n\rangle}$ du produit cartésien $\mathrm{Spec}(B)^{\times_S n}$ sous l'action naturelle du groupe $\mathfrak{S}_n$. Cette construction se recolle : si par exemple $f : X \to S$ est un morphisme quasi-projectif de présentation finie, on définit de même la **puissance symétrique** $n$-ième $(X/S)^{\langle n\rangle}$ de $X/S$ ; sinon, on peut travailler dans la catégorie des espaces algébriques. Comme pour tout quotient par un groupe fini agissant de façon admissible sur un schéma, le morphisme

---

② La notation non standard a pour but d'éviter une possible confusion avec le foncteur des sections globales.

③ À ne pas confondre avec le produit de mélange [Bourbaki A, chap. IV, §5 n°3].

$(X/S)^n \to (X/S)^{\langle n\rangle}$ est fini surjectif, et induit une bijection (topologique) entre points de l'un et $\mathfrak{S}_n$-orbites de points de l'autre, que l'on peut identifier à un multi-ensemble de cardinal $n$ de points de $X$. Pour $X/S$ plat de type fini sur $S$ nœthérien, on déduit de la structure de cogèbre du paragraphe précédent un morphisme fini $\vee_{n_1,n_2} : (X/S)^{\langle n_1\rangle} \times_S (X/S)^{\langle n_2\rangle} \to (X/S)^{\langle n_1+n_2\rangle}$, correspondant sur les points à l'opération d'union de deux multi-ensembles. (Plus généralement, si $\underline{n} = (n_1, \dots, n_r)$ est un $r$-uplet d'entiers naturels, et $n := |\underline{n}| = n_1 + \cdots + n_r$, on note $\vee_{\underline{n}}$ le morphisme $X^{\langle \underline{n}\rangle} := \prod_{i=1}^r X^{\langle n_i\rangle} \to X^{\langle n\rangle}$.)

**A.1.3**. Soit $\Lambda$ un anneau commutatif. Pour toute paire de faisceaux étales de $\Lambda$-modules $\mathscr{F}_1$ sur $(X/S)^{\langle n_1\rangle}$ et $\mathscr{F}_2$ sur $(X/S)^{\langle n_2\rangle}$, on note $\mathscr{F}_1 \vee \mathscr{F}_2$ le faisceau $\vee_{n_1,n_2\star}(\mathscr{F}_1 \boxtimes \mathscr{F}_2)$ sur $(X/S)^{\langle n_1+n_2\rangle}$ ; cette construction s'étend trivialement au cas d'une famille finie $\mathscr{F}_1, \dots, \mathscr{F}_r$ sur $(X/S)^{\langle n_1\rangle}, \dots, (X/S)^{\langle n_r\rangle}$. Dans le cas particulier où $\mathscr{F}_1 = \dots = \mathscr{F}_r$ est un même faisceau $\mathscr{F}$ sur $X$, le faisceau $\vee_{\underline{1}\star}\mathscr{F}^{\boxtimes r}$ ainsi obtenu sur $X^{\langle r\rangle}$ (où $\underline{1}$ désigne le $r$-uplet $(1, \dots, 1)$ et $\vee_{\underline{1}\star} : X^r \to X^{\langle r\rangle}$ le morphisme canonique) est muni d'une action du groupe $\mathfrak{S}_r$. Le faisceau des invariants, $\mathsf{Fix}(\mathfrak{S}_r \circlearrowright \vee_{\underline{1}\star}\mathscr{F}^{\boxtimes r})$, est une **puissance tensorielle symétrique externe** de $\mathscr{F}$; si $X$ est un point, il correspond à la puissance tensorielle symétrique pour les $\Lambda$-modules (A.1.1). (Ce faisceau est noté $\mathsf{TS}^n_{\mathrm{ext}}(\mathscr{F})$ en [SGA 4 XVII, 5.5.7].) Lorsque $\mathscr{F}$ est *plat* sur $\Lambda$, nous le noterons $\mathscr{F}^{\langle n\rangle}$ ; dans le cas général, nous définissons $\mathscr{F}^{\langle n\rangle}$ comme un conoyau (cf. [SGA 4 XVII, 5.5.8]), en faisceautisant l'expression de $\Upsilon^n_\Lambda(M)$ à partir d'une présentation plate.

**A.1.4**. Les foncteurs *non additifs* $\Upsilon^n_\Lambda$ et $\mathscr{F} \mapsto \mathscr{F}^{\langle n\rangle}$ se dérivent, et nous aurons besoin de cette extension, lorsque nous l'appliquerons à des complexes comme $\mathrm{R}\Gamma(X, \mathscr{L}[1])$, concentrés en degrés $[-1, 1]$.

**A.2**. **Cas des courbes.** On suppose dorénavant que $X$ est une courbe quasi-projective lisse sur un corps.

**A.2.1**. Sous cette hypothèse, le morphisme $\mathrm{Div}^n_X = \mathrm{Hilb}^n_X \to (X)^{\langle n\rangle}$, défini en [Bosch, Lütkebohmert et Raynaud 1990, chap. 9, §9.3, p. 253-254] ou [SGA 4 XVII, 6.3.8-9], du schéma des diviseurs effectifs relatifs finis de degré $n$ vers la puissance symétrique, déduit de la propriété universelle du quotient, est un isomorphisme[4] et les morphismes $\vee_{n_1,n_2}$ correspondent à la somme des diviseurs. (Voir aussi [Polishchuk 2003, 16.4], et [Arbarello, Cornalba, Griffiths et Harris 1985, IV.§2] dans le cas analytique). De plus, pour chaque entier $n \geqslant 0$, l'application $\vee_{n,1} : X^{\langle n\rangle} \times X \to X^{\langle n+1\rangle}$ induit pour chaque point fermé $x \in X$ une immersion fermée $X^{\langle n\rangle} \hookrightarrow X^{\langle n+1\rangle}$ correspondant, en terme de diviseurs, à l'application $D \mapsto D + x$.

**A.2.2**. Dans ce paragraphe, nous rappelons pourquoi la puissance symétrique d'une courbe lisse est lisse, en commençant par le cas particulier de la droite affine.
**A.2.2.1**. (Droites affine et projective). Soit $X = \mathbf{A}^1_k = \mathrm{Spec}(k[T])$ la droite affine sur un corps $k$ et, pour chaque $n \geqslant 1$, le produit fibré $X^n := \mathrm{Spec}(k[T_1, \dots, T_n])$. Les $n$ fonctions symétriques élémentaires en les $T_i$ induisent un morphisme $\mathfrak{S}_n$-équivariant $X^n \to \mathbf{A}^n_k := \mathrm{Spec}(k[\sigma_1, \dots, \sigma_n])$, qui se factorise en un isomorphisme $X^{\langle n\rangle} \to \mathbf{A}^n_k$, d'après le théorème des fonctions symétriques. Il n'est pas difficile d'étendre ce résultat à la droite projective : ${\mathbf{P}^1_k}^{\langle n\rangle} = \mathbf{P}^n_k$, l'espace projectif de dimension $n$ sur $k$. (Voir par exemple [Eisenbud et Harris 2016, prop. 10.5].)
En particulier, les schémas ${\mathbf{A}^1_k}^{\langle n\rangle}$ et ${\mathbf{P}^1_k}^{\langle n\rangle}$ sont lisses sur $k$.
**A.2.2.2**. (Critère infinitésimal). Pour montrer que les puissances symétriques de courbes lisses sont lisses, il est possible de se ramener au cas particulier précédent par localisation étale, tout en prenant garde au fait que le foncteur « produit symétrique » ne transforme pas morphisme étale en morphisme étale, cf. A.2.4.2[5]. Nous suivons plutôt une autre approche, qui sera essentielle pour la suite de l'article. D'après ce qui précède, il faut montrer que le schéma $\mathrm{Div}^n_X = \mathrm{Hilb}^n_X$ est lisse, ce qui revient à établir qu'il n'y a pas d'obstruction à relever infinitésimalement un diviseur sur une courbe : c'est vrai par annulation du $H^2$ zariski, cf. p. ex. [Illusie 2005, 8.5.3].

[4] Prendre garde que, sur une base quelconque, la définition de $\mathrm{Div}_{X/S}$ dans [Bosch, Lütkebohmert et Raynaud 1990] ne fait pas d'hypothèse de propreté sur $S$, contrairement à [SGA 4 XVII]. En dimension relative supérieure, le lien entre $\mathrm{Hilb}^n_{X/S}$ et $(X/S)^{\langle n\rangle}$ est plus subtil : cf. p. ex. [Ekedahl et Skjelnes 2014, 7.25].

[5] Contrairement à ce que semble indiquer l'argument [Bosch, Lütkebohmert et Raynaud 1990, p. 255].

**A.2.3.** *Espaces tangents.* Ce même point de vue décrit l'espace tangent — les déformations infinitésimales — comme le $H^0$ à valeurs dans le fibré normal, c'est-à-dire ici :

$$T_D X^{\langle n\rangle} \simeq H^0(X, \mathscr{O}_X(D)/\mathscr{O}_X).$$

où $D$ est un diviseur effectif de degré $n$. (Voir par exemple [Arbarello, Cornalba et Griffiths 2011, chap. IX, §8] et [Polishchuk 2003, prop. 19.2] pour une démonstration plus élémentaire.)

Ce résultat « ponctuel » se globalise de la façon suivante — voir p. ex. [Arbarello, Cornalba, Griffiths et Harris 1985, IV.§2]. Fixons un entier $n \geqslant 1$. Posons $Y := X^{\langle n\rangle}$ et considérons le diviseur universel $\{(x, D) : D - x \geqslant 0\} =: \mathscr{D}_n \hookrightarrow X_Y$ sur la courbe relative $\pi : X_Y \to Y$. Le fibré tangent de $Y$ s'identifie à l'image directe $\pi_\star \mathscr{O}_{\mathscr{D}_n}(\mathscr{D}_n)$ du faisceau $\mathscr{O}_{\mathscr{D}_n}(\mathscr{D}_n) := \mathscr{O}_{X_Y}(\mathscr{D}_n)/\mathscr{O}_{X_Y}$ supporté par le diviseur universel.

**A.2.4.** À des fins de dévissage, nous utiliserons de façon répétée le corollaire suivant, et ses variantes infra, du fait qu'un morphisme quotient par un groupe fini est étale au voisinage d'un point où l'inertie est triviale ([SGA 1 V, 2.2]).

**A.2.4.1. Lemme.** *Soient $X$ une courbe lisse sur un corps algébriquement clos et $D = \sum_{i=1}^r D_i$ un diviseur effectif où les $D_i = n_i x_i$ sont à supports disjoints. Le morphisme $\vee_{\underline{n}} : X^{\langle \underline{n}\rangle} \to X^{\langle n\rangle}$, où $\underline{n} := (n_1, \ldots, n_r)$, est étale au voisinage de $(D_1, \ldots, D_r)$.*

(Ce résultat ramène la lissité de $X^{\langle n\rangle}$ en tout point à celle en les points les plus générés, les $n \cdot x$. Pour une démonstration « algébrique » de la lissité en ces points, voir [Polishchuk 2003, prop. 16.2].)

**A.2.4.2. Lemme.** *Soient $Y \to X$ un morphisme étale entre courbes lisses sur un corps algébriquement clos et $E = \sum_{i=1}^r n_i y_i$ un diviseur effectif de degré $n$ sur $Y$. Le morphisme $Y^{\langle n\rangle} \to X^{\langle n\rangle}$ est étale en $E$ si et seulement si les images $x_i$ des $y_i$ sont distinctes.*

*Démonstration.* Si les images sont distinctes, le lemme précédent nous ramène à montrer que $(Y^{\langle n\rangle}, ny) \to (X^{\langle n\rangle}, nx)$ est étale. Passant aux complétés, cela revient à vérifier que le morphisme, qui est un isomorphisme, $\widehat{\mathscr{O}_{X,x}} \widehat{\otimes}_k \cdots \widehat{\otimes}_k \widehat{\mathscr{O}_{X,x}} \to \widehat{\mathscr{O}_{Y,y}} \widehat{\otimes}_k \cdots \widehat{\otimes}_k \widehat{\mathscr{O}_{Y,y}}$ induit un isomorphisme sur les $\mathfrak{S}_n$-invariants. Réciproquement, en passant aux espaces tangents, il suffit d'observer que si $a > b \geqslant 0$ sont deux entiers et $y \in Y \mapsto x \in X$, le morphisme $H^0(Y, \mathscr{O}_Y(by)/\mathscr{O}_Y) \to H^0(X, \mathscr{O}_X(ax)/\mathscr{O}_X)$ n'est *pas* surjectif, ne serait-ce que pour des raisons de dimension. □

Nous ferons également usage du lemme suivant, de démonstration immédiate (sans hypothèse sur la dimension).

**A.2.4.3. Lemme.** *Soit $X = X_1 \amalg X_2$ le coproduit de deux $k$-schémas quasi-projectifs. Pour tout entier $n \geqslant 1$, les flèches $X_1^{\langle a\rangle} \times X_2^{\langle b\rangle} \to X^{\langle a\rangle} \times X^{\langle a\rangle} \to X^{\langle n\rangle}$ déduites pour chaque décomposition $n = a + b$ des immersions canoniques et de la flèche $\vee_{a,b}$ induisent un isomorphisme*

$$\coprod_{a+b=n} X_1^{\langle a\rangle} \times X_2^{\langle b\rangle} \xrightarrow{\sim} X^{\langle n\rangle}.$$

## A.3. Faisceaux sur une courbe.

**A.3.1.** Commençons par signaler quelques variantes des résultats précédents pour les faisceaux étales.

**A.3.1.1. Lemme.** *Soient $X$ une courbe lisse sur un corps algébriquement clos, $D = \sum_{i=1}^r D_i$ un diviseur effectif où les $D_i = n_i x_i$ sont à supports disjoints. Étale-localement, la paire $(X^{\langle n\rangle}, \mathscr{F}^{\langle n\rangle})$ est isomorphe au voisinage de $D$ à la paire $\big(X^{\langle \underline{n}\rangle} := \prod_{i=1}^r X^{\langle n_i\rangle}, \mathscr{F}^{\langle \underline{n}\rangle} := \boxtimes_{i=1}^r \mathscr{F}^{\langle n_i\rangle}\big)$, au voisinage du point $(D_1, \ldots, D_r)$.*

Passant aux fibres, on en déduit pour $\mathscr{F}$ *plat*, un isomorphisme ([SGA 4 XVII, 5.5.7.2])

$$\mathscr{F}^{\langle n\rangle}{}_y \simeq \bigotimes_{i=1}^r (\mathscr{F}_{y_i})^{\langle n_i\rangle},$$

où $y$ (resp. $y_i$) correspond au diviseur $D$ (resp. $D_i$).

En particulier, on voit que (1) même si $\mathscr{F}$ est lisse, il n'en est pas nécessairement de même de ses puissances tensorielles symétriques externes, sauf en rang 1 ; (2) $\mathscr{F}^{\langle n\rangle}$ est également plat.

**A.3.1.2. Remarques.**

- **(i)** L'observation (1) *supra* vaut également pour les schémas : une puissance symétrique $X^{\langle n\rangle}$ avec $n>1$ d'un schéma lisse $X$ n'est lisse que si $\dim(X)\leqslant 1$, comme il résulte de la pureté et du fait que le lieu de ramification de $X^n\to X^{\langle n\rangle}$ est de codimension $\dim(X)$ dans $X^n$. (Voir aussi [Eisenbud et Harris 2016, prop. 10.6].) Si $X=\mathbb{A}^2_{\mathbb{C}}$ et $n=2$, l'éclatement de la diagonale est une désingularisation, isomorphe au schéma de Hilbert des sous-schémas de longueur 2.
- **(ii)** En rang 1 sur une courbe, on peut vérifier que le faisceau lisse $\mathscr{F}^{\langle n\rangle}$ obtenu correspond à $\mathscr{F}$ via l'isomorphisme $\pi_1(X^{\langle n\rangle})\simeq\pi_1(X)^{\mathrm{ab}}$ [SGA 1 IX, 5.8].
- **(iii)** Notons également que si le rang du faisceau lisse $\mathscr{F}$ sur une courbe $X$ est $>1$, le complexe $\mathscr{F}^{\langle n\rangle}[n]$ est un *faisceau pervers* sur $X^{\langle n\rangle}$ : c'est le prolongement intermédiaire de sa restriction au complémentaire de la diagonale, au-dessus de laquelle c'est un faisceau *lisse* (décalé). (Cf. [Gaitsgory 2007, prop. 5.4].)

**A.3.1.3. Lemme.** *Soient $X=X_1\amalg X_2$ le coproduit de deux $k$-schémas quasi-projectifs et $\mathscr{F}$ un faisceau étale sur $X$ dont on note respectivement $\mathscr{F}_1$ et $\mathscr{F}_2$ les restrictions aux ouverts. Pour chaque décomposition $n=a+b$, la restriction du faisceau $\mathscr{F}^{\langle n\rangle}$ à l'ouvert $X_1^{\langle a\rangle}\times X_2^{\langle b\rangle}$ de $X^{\langle n\rangle}$ est isomorphe à $\mathscr{F}_1^{\langle a\rangle}\boxtimes\mathscr{F}_2^{\langle b\rangle}$.*

**A.3.2.** *Filtrations.*
**A.3.2.1.** Nous utiliserons également de façon répétée l'existence d'un isomorphisme [SGA 4 XVII, 5.5.11.2]

$$(\mathscr{F}\oplus\mathscr{G})^{\langle n\rangle}\simeq\bigoplus_{n_1+n_2=n}\mathscr{F}^{\langle n_1\rangle}\vee\mathscr{G}^{\langle n_2\rangle},$$

lorsque $\mathscr{F}$ et $\mathscr{G}$ sont plats. (Voir aussi [Polishchuk 2003, lemme 19.5] dans le cas cohérent.) Plus généralement, une suite exacte $0\to\mathscr{F}\to\mathscr{E}\to\mathscr{G}\to 0$ de tels faisceaux donne naissance à une filtration à $n+1$ crans de $\mathscr{E}^{\langle n\rangle}$ dont les gradués sont les $\mathscr{F}^{\langle n_1\rangle}\vee\mathscr{G}^{\langle n_2\rangle}$.
**A.3.2.2.** Dans [SGA 4 XVII, 5.5], Deligne fait passer aux objets cosimpliciaux le foncteur (non additif) $\mathscr{F}\longmapsto\mathscr{F}^{\langle n\rangle}$. En utilisant des résolutions plates et l'équivalence de Dold-Kan, il étend ce foncteur à des catégories de complexes. En utilisant les objets bisimpliciaux mixtes de [Illusie 1971-1972, I.4] et en supposant $\Lambda$ noethérien, il est possible de dériver $\mathscr{F}\longmapsto\mathscr{F}^{\langle n\rangle}$ en un foncteur $\mathsf{D}^{\mathrm{b}}_{\mathrm{ctf}}(X,\Lambda)\to\mathsf{D}^{\mathrm{b}}_{\mathrm{ctf}}(X^{\langle n\rangle},\Lambda)$. Si $K\in\mathsf{D}^{\mathrm{b}}_{\mathrm{ctf}}(X,\Lambda)$, on notera encore $K^{\langle n\rangle}\in\mathsf{D}^{\mathrm{b}}_{\mathrm{ctf}}(X^{\langle n\rangle},\Lambda)$ son image par ce foncteur. La décomposition A.3.2.1 vaut encore pour ce foncteur.

Plus généralement, si $K'\to K\to K''\to K'[1]$ est un triangle distingué dans $\mathsf{D}^{\mathrm{b}}_{\mathrm{ctf}}(X,\Lambda)$, alors $K^{\langle n\rangle}$ peut être représenté par un complexe filtré dont les quotients successifs sont les $K'^{\langle n_1\rangle}\vee K''^{\langle n_2\rangle}$ pour $n_1+n_2=n$ (voir [SGA 4 XVII, 5.5.15]). En particulier, le gradué de $K^{\langle n\rangle}$ s'identifie à $(K'\oplus K'')^{\langle n\rangle}$.

## Annexe B. Morphisme d'Abel-Jacobi : rappels

**B.1.** Soit $X$ une courbe projective lisse connexe sur un corps algébriquement clos $k$, dont on note $g$ le genre. Pour tout entier $n\geqslant 1$, notons $\int^n : X^{\langle n\rangle}\to\mathsf{Pic}^n_X$ le morphisme d'Abel-Jacobi envoyant un diviseur effectif $D=x_1+\cdots+x_n$ sur le faisceau $\mathscr{O}(D)$ des fonctions méromorphes $f$ satisfaisant, sur chaque ouvert, $(f)+D\geqslant 0$. D'après un théorème d'Abel, la fibre du morphisme, propre, $\int^n$ au-dessus de la classe $[\mathscr{L}]$ d'un faisceau inversible $\mathscr{L}$ est en bijection est le système linéaire complet $\mathbb{P}(H^0_{\mathrm{coh}}(X,\mathscr{L})^\vee)$, lorsque $h^0_{\mathrm{coh}}(X,\mathscr{L})\neq 0$. Notons $r(\mathscr{L}):=h^0_{\mathrm{coh}}(X,\mathscr{L})-1$ la dimension de cet espace projectif. On a l'encadrement trivial

$$\dim X^{\langle n\rangle}-\dim\mathsf{Pic}^n_X=n-g\leqslant r\leqslant n=\dim X^{\langle n\rangle},$$

où $g$ est le genre de la courbe $X$. Le théorème de Riemann-Roch apporte des informations supplémentaires sur ces fibres : la caractéristique d'Euler-Poincaré cohérente $\chi_{\mathrm{coh}}(X,\mathscr{L}):=h^0_{\mathrm{coh}}(X,\mathscr{L})-h^1_{\mathrm{coh}}(X,\mathscr{L})$ est égale à $\deg(\mathscr{L})+\chi(\mathscr{O}_X)$, ou encore $\deg(\mathscr{L})-g+1$. En particulier, $h^0_{\mathrm{coh}}(X,\mathscr{L})>0$ si $\deg(\mathscr{L})\geqslant g$. D'autre part, il résulte de la dualité de Poincaré, dans le cas cohérent, que $h^1_{\mathrm{coh}}(X,\mathscr{L})=h^0_{\mathrm{coh}}(X,\omega_X\otimes\mathscr{L}^\vee)$.

L'acyclicité, pour la topologie étale, d'une paire (morphisme, faisceau) étant souvent une incarnation cohomologique de la lissité des deux constituants de la paire, il est naturel de commencer par s'intéresser à la géométrie des morphismes d'Abel-Jacobi, à la fois en le degré qui nous intéresse mais aussi en degrés inférieurs : les puissances symétriques plus petites, que l'on peut voir (non canoniquement) comme des sous-schémas fermés de la puissance symétrique originale, étant souvent des strates de lissité de sous-quotients du faisceau puissance symétrique.

**B.2.** Pour $n$ grand, le morphisme d'Abel-Jacobi est lisse : plus précisément, lorsque $n > 2g-2 = \deg(K) = \chi(X)$, on a $r = n-g$ et $\int^n$ est surjectif et lisse. D'après un théorème de Jacobi, la surjectivité persiste pour $n \geqslant g$ et on a donc *en général* $r = n-g$. Si $n < g$, on a par contre *en général* $r = 0$ : le morphisme d'Abel-Jacobi est génériquement fini sur son image. (C'est également vrai pour $n = g$.) La lissité pour $n > 2g-2$ peut être énoncée sous la forme géométrique, plus précise, suivante : la puissance symétrique $X^{\langle n \rangle}$ est le fibré projectif sur $\mathrm{Pic}^n_X$ d'un faisceau inversible de degré $n$.

**B.3.** Soit $D$ un diviseur effectif de degré $n$ sur $X$. La description faite en A.2.3 de l'espace tangent $T_D X^{\langle n \rangle}$ permet d'identifier la différentielle $T_D \int^n$ au cobord $H^0_{\mathrm{coh}}(X, \mathscr{O}_D(D)) \to H^1_{\mathrm{coh}}(X, \mathscr{O}_X)$ associé à la suite exacte $0 \to \mathscr{O}_X \to \mathscr{O}_X(D) \to \mathscr{O}_D(D) \to 0$ ; voir [Polishchuk 2003, prop. 19.2]. On utilise ici le fait que $T_{[\mathscr{O}(D)]} \mathrm{Pic}^n_X$ est isomorphe à Lie $\mathrm{Pic}^0_X = H^1_{\mathrm{coh}}(X, \mathscr{O}_X)$. Par dualité, l'image de $T_D \int^n$ s'identifie donc au noyau de $H^0_{\mathrm{coh}}(X, \Omega^1)^\vee \to H^0_{\mathrm{coh}}(X, \Omega^1(-D))^\vee$, c'est-à-dire à l'orthogonal de $H^0_{\mathrm{coh}}(X, \Omega^1(-D)) \subseteq H^0_{\mathrm{coh}}(X, \Omega^1)$ dans $H^0_{\mathrm{coh}}(X, \Omega^1)^\vee \simeq T_{[\mathscr{O}(D)]} \mathrm{Pic}^n_X$.

# Bibliographie

## Sigles

### Éléments de mathématique

### Séminaires de géométrie algébrique

## Autres références